# Cartan Angular Invariant and Deformations in Rank One Symmetric Spaces

Boris Apanasov and Inkang Kim

**Abstract**

We develop and study quaternionic and octonionic analogies of Cartan angular and Toledo invariants that are well known in the complex hyperbolic space. Using such invariants we study quasifuchsian deformations (including bendings) of quaternionic and octonionic hyperbolic manifolds.

## 1 Introduction

Many authors studied deformations of locally symmetric spaces of rank one, that is spaces modeled on $\mathbb{F}$-hyperbolic spaces $H_\mathbb{F}^n$ over numbers $\mathbb{F}$ that are either real $\mathbb{R}$, complex $\mathbb{C}$, quaternions $\mathbb{H}$, or octonions $\mathbb{O}$. Besides Mostow [23] rigidity theorem, Corlette [12] and Gromov–Schoen [19] extended the Margulis [22] superrigidity theorem for lattices in semisimple Lie groups of real rank at least two to the case of real rank one which corresponds to automorphisms groups $O(n,1)$, $U(n,1)$, $Sp(n,1)$ and $F_4^{-20}$ of real, complex and quaternionic hyperbolic spaces and the hyperbolic Cayley plane. Namely they proved that any lattice $\Gamma$ in $Sp(n,1)$, $n \geq 2$, or $F_4^{-20}$ is superrigid over archimedian fields and in $p$-adic case (which also implies its arithmeticity). For a geometric sense of such superrigidity for quaternionic manifolds, see Kamishima [20] and Apanasov-Kim [5].

In the remaining cases of real and complex hyperbolic spaces we have many non-arithmetic lattices and there are a number of constructions which show that superrigidity does not hold here either, see Apanasov [1] for real

---

[1]1991 *Mathematics Subject Classification.* 51M10, 57S25.

hyperbolic lattices and Apanasov [3] for bending deformations in the complex hyperbolic space.

On the other hand, W.Goldman [15] discovered a new type of rigidity for embeddings of a closed surface group $\pi_1 S = \Gamma \to PU(2,1)$ nearby its inclusion $\Gamma \subset PU(1,1) \subset PU(2,1)$. A reasoning behind this result was the Toledo invariant for complex hyperbolic surfaces $M$ homotopy equivalent to a closed Riemann surface $S$, $f: S \to M$, defined as the integral over $S$ of the pull back $f_*\omega$ of the Kähler form $\omega$ on $M$. Due to Toledo [25], $f(S) \subset M$ is complex geodesic if that invariant equals $\pm 2\pi\chi(S)$. This is the special case of the Goldman-Millson's [18] conjecture that any representation $\rho \in \text{Hom}(\Gamma, SU(n+k,1))$ of a uniform lattice $\Gamma \subset SU(n,1)$ is a faithful discrete representation stabilizing a totally geodesic $H_{\mathbb{C}}^n$ in $H_{\mathbb{C}}^{n+k}$, provided that $\int_{H_{\mathbb{C}}^n/\Gamma} \rho^*\omega^n = \text{Vol}(H_{\mathbb{C}}^n/\Gamma)$ where $\omega$ is a Kähler form. The conjecture has been settled by Corlette in [11], and due to a different proof by [18], the above local rigidity appears to be true in all dimensions:

*Let $\Gamma \subset PU(n-1,1) \subset PU(n,1)$, $n > 1$, be a uniform lattice in $PU(n-1,1)$, $n \geq 2$. Then for any its representation $\rho: \Gamma \to PU(n,1)$ nearby the inclusion, the group $\rho(\Gamma)$ preserves a complex totally geodesic (n-1)-subspace where its action is conjugate to that of $\Gamma$.*

Moreover in dimensions $n > 2$, this local rigidity of complex hyperbolic $n$-manifolds $M$ homotopy equivalent to their closed complex totally geodesic hypersurfaces is even global (at least in the connected component of representation variety [11]). We note however that both conditions of this rigidity, the action of a lattice in an *analytic* subspace and its co-compactness are essential. Firstly, there are quasifuchsian deformations of non-uniform lattices in complex hyperbolic geometry which appear to be quasiconformally instable, see Apanasov [4], Falbel-Koseleff [14]. And secondly, non-rigidity of uniform lattices acting in totally real totally geodesic subspaces in $H_{\mathbb{C}}^2$ follows from Apanasov's [3] construction of bending deformations in the complex hyperbolic space.

The goal of this paper is to address these deformation questions in the case of quaternionic and octonionic hyperbolic structures. Namely we construct nontrivial deformations of quaternionic and octonionic hyperbolic manifolds of infinite volume. These deformations correspond to quasifuchsian deformations of (Fuchsian) inclusions of the fundamental groups of such manifolds into $Sp(n,1)$, $n \geq 2$, or $F_4^{-20}$, in particular bending deformations of the inclusions. We will address new rigidity effects (see [5]) in another paper.

To distinguish different quasifuchsian representations and to prove their non-triviality, we use some new invariants for triples of points in the closure $\overline{H_\mathbb{F}^n} = H_\mathbb{F}^n \cup \partial H_\mathbb{F}^n$ of either quaternionic hyperbolic n-space or the Cayley hyperbolic plane. We develop these invariants by using corresponding Hermitian products and Kähler (quaternionic or octonionic) four- and eight-forms. So these invariants are some generalizations of well known in the complex hyperbolic geometry Cartan angular invariant and Toledo invariant, see Goldman [16].

**Acknowledgements.** This work is partially supported by Korea Advanced Institute of Science and Technology, the University of Oklahoma and Université de Paris-Sud. All their partial support is greatly acknowledged. We have enjoyed fruithful discussions with many experts including Werner Ballmann, Vicente Cortés, Misha Gromov, Ursula Hamenstädt, Yoshi Kamishima, Pierre Pansu, Frederic Paulin and Larry Siebenmann. We would like to use this opportunity to thank all of them.

## 2 Preliminaries

### 2.1 Symmetric spaces of rank one

The symmetric spaces of $\mathbb{R}$-rank one of non-compact type are the hyperbolic spaces $H_\mathbb{F}^n$, where $\mathbb{F}$ is either the real numbers $\mathbb{R}$, or the complex numbers $\mathbb{C}$, or the quaternions $\mathbb{H}$, or the Cayley numbers $\mathbb{O}$; in last case $n = 2$. They are respectively called as real, complex, quaternionic and octonionic hyperbolic spaces (the latter one $H_\mathbb{O}^2$ is also known as the Cayley hyperbolic plane). Algebraically these spaces can be described as the corresponding quotients: $SO(n,1)/SO(n)$, $SU(n,1)/SU(n)$, $Sp(n,1)/Sp(n)$ and $F_4^{-20}/Spin(9)$ where the latter group $F_4^{-20}$ of automorphisms of the Cayley plane $H_\mathbb{O}^2$ is the real form of $F_4$ of rank one. We normalize the metric so the (negative) sectional curvature of $H_\mathbb{F}^n$ is bounded from below by $-1$.

Following Mostow [23] and using the standard involution (conjugation) in $\mathbb{F}$, $z \to \bar{z}$, one can define projective models of the hyperbolic spaces $H_\mathbb{F}^n$ as the set of negative lines in the Hermitian vector space $\mathbb{F}^{n,1}$, with Hermitian structure given by the indefinite $(n,1)$-form

$$\langle\langle z, w \rangle\rangle = z_1 \overline{w}_1 + \cdots + z_n \overline{w}_n - z_{n+1} \overline{w}_{n+1}.$$

Here, taking non-homogeneous coordinates, one can obtain unit ball models

(in the unit ball $B^n_{\mathbb{F}}(0,1) \subset \mathbb{F}^n$) for the first three spaces (Since the multiplication by quaternions is not commutative, we specify that we use "left" vector space $\mathbb{H}^{n,1}$ where the multiplication by quaternion numbers is on the left). However, it does not work for the Cayley plane since $\mathbb{O}$ is non-associative, and one should use a Jordan algebra of $3 \times 3$ Hermitian matrices with entries from $\mathbb{O}$ whose group of automorphisms is $F_4$, see [23].

Another models of $H^n_{\mathbb{F}}$ use the so called horospherical coordinates [16, 7] based on foliations of $H^n_{\mathbb{F}}$ by horospheres centered at a fixed point $\infty$ at infinity $\partial H^n_{\mathbb{F}}$ which is homeomorphic to $(n \dim_{\mathbb{R}} \mathbb{F} - 1)$-dimensional sphere. Such a horosphere can be identified with the nilpotent group $N$ in the Iwasawa decomposition $KAN$ of the automorphism group of $H^n_{\mathbb{F}}$. The nilpotent group $N$ can be identified with the product $\mathbb{F}^{n-1} \times \operatorname{Im} \mathbb{F}$ (see [23]) equipped with the operations:

$$(\xi, v) \cdot (\xi', v') = (\xi + \xi', v + v' + 2 \operatorname{Im} \langle \xi, \xi' \rangle) \quad \text{and} \quad (\xi, v)^{-1} = (-\xi, -v),$$

where $\langle , \rangle$ is the standard Hermitian product in $\mathbb{F}^{n-1}$, $\langle z, w \rangle = \sum z_i \overline{w_i}$. The group $N$ is a 2-step nilpotent Carnot group with center $\{0\} \times \operatorname{Im} \mathbb{F} \subset \mathbb{F}^{n-1} \times \operatorname{Im} \mathbb{F}$, and acts on itself by the left translations $T_h(g) = h \cdot g$, $h, g \in N$.

Now we may identify

$$H^n_{\mathbb{F}} \cup \partial H^n_{\mathbb{F}} \backslash \{\infty\} \longrightarrow N \times [0, \infty) = \mathbb{F}^{n-1} \times \operatorname{Im} \mathbb{F} \times [0, \infty),$$

and call this identification the *"upper half-space model"* for $H^n_{\mathbb{F}}$ with the natural horospherical coordinates $(\xi, v, u)$. In these coordinates, the above left action of $N$ on itself extends to an isometric action (Carnot translations) on the $\mathbb{F}$-hyperbolic space in the following form:

$$T_{(\xi_0, v_0)} : (\xi, v, u) \longmapsto (\xi_0 + \xi, v_0 + v + 2 \operatorname{Im} \langle \xi_0, \xi \rangle, u),$$

where $(\xi, v, u) \in \mathbb{F}^{n-1} \times \operatorname{Im} \mathbb{F} \times [0, \infty)$.

There are a natural norm and an induced by this norm distance on the Carnot group $N = \mathbb{F}^{n-1} \times \operatorname{Im} \mathbb{F}$, which are known in the case of the Heisenberg group (when $\mathbb{F} = \mathbb{C}$) as the Cygan's norm and distance. Using horospherical coordinates, they can be extended to a norm on $H^n_{\mathbb{F}}$, see [21, 7]:

$$|(\xi, v, u)|_c = |(|\xi|^2 + u - v)|^{1/2}, \tag{1}$$

where $|.|$ is the norm in $\mathbb{F}$, and to a metric $\rho_c$ (still called the Cygan metric) on the upper half-space model $\mathbb{F}^{n-1} \times \operatorname{Im} \mathbb{F} \times (0, \infty)$ of $H^n_{\mathbb{F}}$:

$$\rho_c\big((\xi,v,u),(\xi',v',u')\big) = \big|\,|\xi-\xi'|^2 + |u-u'| - (v-v'+2\operatorname{Im}\langle\xi,\xi'\rangle)\,\big|^{\frac{1}{2}}. \quad (2)$$

It follows directly from the definition that Carnot translations and rotations are isometries with respect to the Cygan metric $\rho_c$. Moreover, the restrictions of this metric to different horospheres centered at $\infty$ are the same, so Cygan metric plays the same role as Euclidean metric does on the upper half-space model for the real hyperbolic space $\mathbb{H}^n$.

There is a projection from $N \cup \infty = \partial(\mathbb{F}^{n-1} \times \operatorname{Im}\mathbb{F} \times (0,\infty))$ to the boundary of the unit ball in $\mathbb{F}^n$, which is the extension to infinity of an isometry between two models of $H_\mathbb{F}^n$ and has the form (see [21]):

$$[z,t] \to \left[2\frac{1+|z|^2+t}{|1+|z|^2+t|^2}z,\ \frac{1+|z|^2+t}{|1+|z|^2+t|^2}(1-|z|^2+t)\right] \in S(0,1) \subset \mathbb{F}^n, \quad (3)$$

where $z \in \mathbb{F}^{n-1}$ and $t \in \operatorname{Im}\mathbb{F}$

Note that two points, the origin $[0,0] \in N$ and $\infty$ correspond to $(0,1)$ and $(0,-1)$ respectively where $(0,\pm 1) = (0,\cdots,0,\pm 1) \in \mathbb{F}^n$. A real hyperbolic two plane in the unit ball model of $H_\mathbb{F}^n$ whose boundary at infinity contains the points $(0,1)$ and $(0,-1)$ corresponds to a 2-plane in $\mathbb{F}^{n-1} \times \operatorname{Im}\mathbb{F} \times (0,\infty)$ whose boundary line in $N$ passes through the origin and some point $[z,0]$. The reason is as follows. The real circle $\{[\mathbb{R},0]\} \cup \{\infty\}$ is obviously the boundary of a real hyperbolic two plane. Now the group $K_0 A$ fixing the origin of $N$ and $\infty$ is $U(n-1) \times \mathbb{R} \times U(1)$ in $H_\mathbb{C}^n$, $Sp(n-1) \times \mathbb{R} \times Sp(1)$ in $H_\mathbb{H}^n$, and $Spin(7) \times \mathbb{R}$ in $H_\mathbb{O}^2$. Then for any $[z,0]$ with $z \in \mathbb{F}^{n-1}$, $|z|=1$, there is a transformation $M \in K_0$ which maps $(1,0,...,0) \in \mathbb{F}^{n-1}$ to $z$. So $M$ maps $[1,0]$ to $[z,0]$, which implies that it maps $\{[\mathbb{R},0]\}$ to the line in $N$ that passes through the origin and $[z,0]$.

We can define an orthogonal projection to an $\mathbb{F}$-subspace $L \subset H_\mathbb{F}^n$, $\Pi\colon H_\mathbb{F}^n \to L$. Indeed, since the metric in $H_\mathbb{F}^n$ is strictly negatively curved, its distance function is strictly convex. Hence for any $x \in H_\mathbb{F}^n$ there is a unique point $\Pi(x)$ in $L$ closest to $x$. For any points $u \in H_\mathbb{F}^n$ and $v \in L$, the (real) geodesics from $\Pi(u)$ to $u$ and to $v$ are orthogonal and span a totally real 2-plane, see Lemma 2.2.

The group of automorphisms of $H_\mathbb{H}^n$ is $PSp(n,1)$. The stabilizer $K$ of the origin is $Sp(1) \times Sp(n)$ which can be described in the matrix form as:

$$\begin{bmatrix} M & 0 \\ 0 & \nu \end{bmatrix}$$

where $M \in Sp(n)$ and $\nu \in Sp(1)$. Note the matrix acts on the right, and the projectivization is given by multiplication on the left. So in the ball model, the action is:
$$z \to \nu^{-1} z M .$$

The stabilizer of a real geodesic connecting two points $(0, 1)$ and $(0, -1)$ is $MA = Sp(1) \times Sp(n-1) \times \mathbb{R}$. This action can be described in the matrix form as:
$$\begin{bmatrix} M & 0 & 0 \\ 0 & \nu \cosh r & \nu \sinh r \\ 0 & \nu \sinh r & \nu \cosh r \end{bmatrix}$$
where $M \in Sp(n-1)$, $\nu \in Sp(1)$ and $r \in \mathbb{R}$. Specially $Sp(1)$ acts as
$$\begin{bmatrix} I & 0 & 0 \\ 0 & \nu & 0 \\ 0 & 0 & \nu \end{bmatrix}$$

If $(0, \mathbb{H})$ is the $\mathbb{H}$-line containing the real geodesic joining $(0, 1)$ and $(0, -1)$, the action of $Sp(1)$ on this $\mathbb{H}$-line is:
$$(0, \mathbb{H}) \to \nu^{-1}(0, \mathbb{H}\nu, \nu) = (0, \nu^{-1}\mathbb{H}\nu).$$

But in general, $\nu \in Sp(1)$ maps $(z, z_n)$ to $(\nu^{-1} z, \nu^{-1} z_n \nu)$.

For a geodesic triangle $\triangle$ in the $\mathbb{H}$-line whose one edge is contained in the real axis, its orbit under $Sp(1)$ is a four dimensional object in the $\mathbb{H}$-line. We will use this four dimensional object to define the Toledo invariant in §4.

A Cayley number $z \in \mathbb{O}$ is a pair of quaternions, $z = (q_1, q_2)$, and the multiplication in $\mathbb{O}$ is given by
$$(q_1, q_2)(p_1, p_2) = (q_1 p_1 - \bar{p}_2 q_2, p_2 q_1 + q_2 \bar{p}_1) .$$

The standard involution (conjugation) in $\mathbb{O}$ is defined by $\overline{(q_1, q_2)} = (\bar{q}_1, -q_2)$, so for $z = (q_1, q_2) \in \mathbb{H} \times \mathbb{H} = \mathbb{O}$, we have $\operatorname{Im} z = (\operatorname{Im} q_1, q_2)$ and $\operatorname{Re} z = \operatorname{Re} q_1$. Then Cayley numbers satisfy the usual properties like: $x\bar{x} = |x|^2$, $|xy| = |x||y|$, $x^{-1} = \bar{x}/|x|^2$, $\overline{xy} = \bar{y}\bar{x}$. Even though Cayley numbers are not commutative, nor associative, by Artin's lemma a subalgebra generated by two elements is associative. Cayley hyperbolic plane is made out of an exceptional Jordan algebra of $3 \times 3$ Hermitian matrices with entries from $\mathbb{O}$ whose group of automorphisms is $F_4$, see [23]. The group of automorphisms of the Cayley plane $H_\mathbb{O}^2$ is $F_4^{-20}$, the real form of $F_4$ of rank one. The stabilizer in $F_4^{-20}$ of the origin $(0,0) \in B_{\mathbb{O}}^2(0,1) = H_{\mathbb{O}}^2$ is $Spin(9)$ operating on $\mathbb{O}^2 = \mathbb{R}^{16}$ via the spinor representation. If $L_1 = \mathbb{O} \times 0$ and $L_2 = 0 \times \mathbb{O}$ denote the coordinate $\mathbb{O}$-axes, then the stabilizer of $L_1$ acts on $L_1$ as $SO(8)$ via the even $\frac{1}{2}$-spin representation, and on $L_2$ as odd $\frac{1}{2}$-spin representation. The stabilizer of the real line through $(0,0)$ and $(1,0)$ is $Spin(7)$.

## 2.2 Bisectors and Dirichlet polyhedra

It is well-known that there is no totally geodesic (real) hypersurfaces in rank one symmetric spaces $H_{\mathbb{F}}^n$ different from real hyperbolic spaces $H_{\mathbb{R}}^n$, that is for $\mathbb{F} \neq \mathbb{R}$. A reasonable substitute are the bisectors introduced by G.Mostow in the complex hyperbolic geometry, where he pointed out that such bisectors are as close to being totally geodesic as possible. They are minimal hypersurfaces of cohomogenity $(\dim_{\mathbb{R}} \mathbb{F} - 1)$, all equivalent under isometries of $H_{\mathbb{F}}^n$, and have a natural decomposition into totally geodesic $\mathbb{F}$-hypersurfaces which we shall describe below.

For two distinct points $z_1, z_2 \in H_{\mathbb{F}}^n$, the *bisector equidistant from $z_1$ and $z_2$* is defined as:

$$\mathfrak{S}(z_1, z_2)) = \{z \in H_{\mathbb{F}}^n : d(z, z_1) = d(z, z_2)\}. \tag{4}$$

The boundary of a bisector at infinity $\partial H_{\mathbb{F}}^n$ is called *a spinal sphere* in $\partial H_{\mathbb{F}}^n = N \cup \{\infty\}$. Following the terminology and arguments in the complex hyperbolic spaces due to Mostow and Goldman [16], we shall give "Mostow's slice decomposition" of bisectors. For given $z_1, z_2 \in H_{\mathbb{F}}^n$ and their bisector $\mathfrak{S}(z_1, z_2)$, let $\Sigma \subset H_{\mathbb{F}}^n$ be the $\mathbb{F}$-line spanned by those points and called the $\mathbb{F}$-spine of $\mathfrak{S}$ (with respect to $z_1$ and $z_2$). Also we define the spine of $\mathfrak{S}$ (with respect to $z_1$ and $z_2$) as

$$\sigma(z_1, z_2) = \mathfrak{S}(z_1, z_2) \cap \Sigma = \{w \in \Sigma : d(z_1, w) = d(z_2, w)\}, \tag{5}$$

which is the orthogonal bisector of the geodesic segment joining $z_1$ and $z_2$ in the sense of the real hyperbolic geometry induced by $H_{\mathbb{F}}^n$ on its $\mathbb{F}$-line $\Sigma$.

**Theorem 2.1** (Mostow). *Let $\mathfrak{S} = \mathfrak{S}(z_1, z_2)$, its $\mathbb{F}$-spine $\Sigma$ and the spine $\sigma$ be as above, and $\Pi_{\Sigma}: H_{\mathbb{F}}^n \to \Sigma$ be the orthogonal projection onto $\Sigma$. Then*

$$\mathfrak{S} = \Pi_{\Sigma}^{-1}(\sigma) = \bigcup_{s \in \sigma} \Pi_{\Sigma}^{-1}(s).$$

**Proof:** Due to the strict convexity of the distance function in $H_\mathbb{F}^n$, for any $z \in H_\mathbb{F}^n \setminus \Sigma$ there is a unique projection point $\Pi(z)$ in $\Sigma$ closest to $z$. As the standard fact of Riemannian geometry, we have that for any point $s \in \Sigma$, the (real) geodesics from $\Pi(z)$ to $z$ and to $s$ are orthogonal. Indeed, if those geodesics are not orthogonal, we can decrease the distance from $z$ to $\Pi(z)$ along the geodesic from $\Pi(z)$ to $s$ until the angle is $\pi/2$. This contradicts to the fact that $\Pi(z)$ is the closest point in $\Sigma$ to $z$. Now, applying the Pythagorean theorem and repeating arguments in $H_\mathbb{C}^n$ (see [16], Lemma 3.2.13), we have:

**Lemma 2.2** *Let $\Sigma \subset H_\mathbb{F}^n$ be a $\mathbb{F}$-subspace with orthogonal projection $\Pi$. Then for all $z \in H_\mathbb{F}^n \setminus \Sigma$ and $s \in \Sigma$, the real geodesics from $\Pi(z)$ to $z$ and to $s$ are orthogonal and span a totally real 2-plane. Furthermore*

$$\cosh\left(\frac{d(z,s)}{2}\right) = \cosh\left(\frac{d(z,\Pi(z))}{2}\right) \cosh\left(\frac{d(\Pi(z),s)}{2}\right).$$

*Conversely, if $z, w, s$ span a right triangle on a totally real 2-plane in $H_\mathbb{F}^n$ then the orthogonal projection to the $\mathbb{F}$-line spanned by $w$ and $s$ maps $z$ to $w$, and the orthogonal projection to the $\mathbb{F}$-line spanned by $z$ and $w$ maps $s$ to $w$.*

Now, applying this Lemma to both our points $z_1$ and $z_2$, we have:

$$\begin{aligned} z \in \mathfrak{S}(z_1, z_2) &\iff d(z, z_1) = d(z, z_2) \\ &\iff d(\Pi_\Sigma(z), z_1) = d(\Pi_\Sigma(z), z_2) \\ &\iff \Pi_\Sigma(z) \in \sigma(z_1, z_2) \,, \end{aligned}$$

which finishes the proof of theorem. ■

Following Mostow, we call the $\mathbb{F}$-hyperplanes $\Pi_\Sigma^{-1}(s)$ (the fibers of the real analytic fibration given by the orthogonal projection $\Pi_\Sigma : H_\mathbb{F}^n \to \Sigma$) by *slices* of $\mathfrak{S}$, their boundary spheres at infinity by $\mathbb{F}$-*chains*, and the boundary of the spine $\sigma \subset \mathfrak{S}$ by the *vertex sphere* (whose dimension is $\dim \mathbb{F} - 2$).

Then, since orthogonal projection $\Pi_\Sigma : H_\mathbb{F}^n \to \Sigma$ is a real analytic fibration, the above theorem immediately implies (cf. [16]):

**Corollary 2.3** *A bisector $\mathfrak{S}$ is a real analytic hypersurface in $H_\mathbb{F}^n$ fibered by its slices. The spinal sphere $\partial \mathfrak{S}$ is a real analytic real hypersurface in $\partial H_\mathbb{F}^n$ fibered by $\mathbb{F}$-chains.*

In particular taking $\partial H_\mathbb{F}^n$ as $(\mathbb{F}^{n-1}\times\operatorname{Im}\mathbb{F})\cup\{\infty\}$, we may represent infinite spinal spheres passing through the origin as the products $\mathbb{F}^{n-1} \times \mathbb{R}^{\dim\mathbb{F}-2}$ where $\mathbb{R}^{\dim\mathbb{F}-2}$ is a linear subspace in $\operatorname{Im}\mathbb{F}$ (its union with $\{\infty\}$ is the vertex sphere). Here we have the natural fibration of the factor $\mathbb{F}^{n-1}$ (and hence of the whole infinite spinal sphere) by $\mathbb{F}$-chains represented by concentric spheres in $\mathbb{F}^{n-1}$ centered at the origin.

By using bisectors we can construct fundamental polyhedra for discrete isometry group actions in a rank one symmetric space $H_\mathbb{F}^n$, the so called Dirichlet polyhedra. Namely the Dirichlet polyhedron of a discrete group $\Gamma \subset \operatorname{Isom} H_\mathbb{F}^n$ centered at a point $z \in H_\mathbb{F}^n$ that is not fixed by non-trivial elements of $\Gamma$ is given by

$$D_z(\Gamma) = \{x \in H_\mathbb{F}^n : d(x,z) < d(x,\gamma(z)),\ \gamma \in \Gamma\backslash\{\mathrm{id}\}\}, \qquad (6)$$

and is bounded by bisectors $\mathfrak{S}(z,\gamma(z))$ equidistant from $z$ and $\gamma(z)$, see [13].

For (real) Fuchsian groups $\Gamma \subset PO(n,1) \subset \operatorname{Isom} H_\mathbb{F}^n$ preserving a totally real n-space $H_\mathbb{R}^n \subset H_\mathbb{F}^n$, behavior of Dirichlet bisectors $\mathfrak{S}(z,\gamma(z)) \subset H_\mathbb{F}^n$, $z \in H_\mathbb{R}^n$, $\gamma \in \Gamma$, is determined by their intersection with the real ($\Gamma$-invariant) n-space $H_\mathbb{R}^n$. In other words, the Dirichlet polyhedron $D_z(\Gamma) \subset H_\mathbb{F}^n$ centered at a point $z \in H_\mathbb{R}^n$ "has the same combinatorics" as its restriction to $H_\mathbb{R}^n$. Namely (assuming $n = 2$), for any two points $z, w \in H_\mathbb{F}^2$, let us denote by $\mathfrak{S}^+(z,w)$ the equidistant half-space

$$\mathfrak{S}^+(z,w) = \{y \in H_\mathbb{F}^2 : d(y,z) < d(y,w)\}. \qquad (7)$$

Then repeating word-by-word W.Goldman's [16] arguments in the complex hyperbolic space, we have

**Theorem 2.4** *For any three points $w, z_1, z_2 \in H_\mathbb{R}^2$, the half-spaces $\mathfrak{S}^+(z_1,w)$ and $\mathfrak{S}^+(z_2,w)$ are disjoint if and only if their intersections $\mathfrak{S}^+(z_i,w) \cap H_\mathbb{R}^2$ are disjoint.*

We finish this section by noticing that in general, for $\mathbb{F} \neq \mathbb{R}$, the image of a geodesic under the orthogonal projection is not geodesic. The simplest counterexample is as follows. Let $\mathfrak{S}(A,B) \subset H_\mathbb{F}^n$ be the bisector of two points $A$ and $B$. Then there are a unique $\mathbb{F}$-line $\Sigma$ containing $A$ and $B$, a unique real geodesic connecting $A$ and $B$ (on the $\mathbb{F}$-line $\Sigma$) and a unique totally geodesic $(k-1)$-plane $\sigma \cong H_\mathbb{R}^{k-1} \subset \Sigma$ bisecting $A$ and $B$. Here the metric of the $\mathbb{F}$-line $\Sigma$ induced by $H_\mathbb{F}^n$ coincides with the Poincaré metric of the real hyperbolic $k$-space $H_\mathbb{R}^k$ with $k = \dim_\mathbb{R} \mathbb{F}$. Each fiber of a point $x \in \sigma$ under the orthogonal projection $\Pi \colon H_\mathbb{F}^n \to \Sigma$ is a totally geodesic $\mathbb{F}$-submanifold of real codimension $k$ which meet the $\mathbb{F}$-line $\Sigma$ orthogonally at the point $x$. This fibration along the real geodesic $(k-1)$-plane $\sigma$ is the bisector $\mathfrak{S}(A, B)$ of two points $A$ and $B$. Since it is not totally geodesic we can choose two points $p, q \in \mathfrak{S}(A, B)$ over $x, y \in \sigma$ respectively so that the geodesic $\gamma$ joining them does not lie on $\mathfrak{S}(A, B)$. Then $\Pi(\gamma) \subset \Sigma$ is a curve between $x$ and $y$ which does not lie in the real geodesic plane $\sigma$, so it is not a geodesic. Finally we note that if the projection of a geodesic is again a geodesic, it is not hard to see that these two geodesics span a totally real two plane.

## 3 Quaternionic Cartan angular invariant

In the complex hyperbolic geometry there is a well known invariant detecting whether a triple of points lies on a $\mathbb{C}$-chain or an $\mathbb{R}$-circle, the Cartan angular invariant, see [16]. Here we would like to extend these ideas to define its analogue in the quaternionic hyperbolic space $H_\mathbb{H}^n$.

Let $x = (x_1, x_2, x_3) \in (H_\mathbb{H}^n \cup \partial H_\mathbb{H}^n) \times (H_\mathbb{H}^n \cup \partial H_\mathbb{H}^n) \times (H_\mathbb{H}^n \cup \partial H_\mathbb{H}^n)$ be a triple of distinct points and $\tilde{x}_i \in H_\mathbb{H}^{n,1}$ some lifts of these points $x_i$.

**Definition 3.1** *The quaternionic Cartan angular invariant of a triple $x$, $0 \leq \mathbb{A}_\mathbb{H}(x) \leq \pi/2$, is the angle between the first coordinate line $\mathbb{R}e_1 = (\mathbb{R}, 0, 0, 0) \subset \mathbb{R}^4$ and the Hermitian triple product*

$$\langle \tilde{x}_1, \tilde{x}_2, \tilde{x}_3 \rangle = \langle \tilde{x}_1, \tilde{x}_2 \rangle \langle \tilde{x}_2, \tilde{x}_3 \rangle \langle \tilde{x}_3, \tilde{x}_1 \rangle \in \mathbb{H},$$

*where we identify $\mathbb{H}$ and $\mathbb{R}^4$.*

If we choose another lifts $\nu_i \tilde{x}_i$ of the points $x_i$, $i = 1, 2, 3$, then

$$\langle \nu_1 \tilde{x}_1, \nu_2 \tilde{x}_2, \nu_3 \tilde{x}_3 \rangle = |\nu_2|^2 |\nu_3|^2 \nu_1 \langle \tilde{x}_1, \tilde{x}_2 \rangle \langle \tilde{x}_2, \tilde{x}_3 \rangle \langle \tilde{x}_3, \tilde{x}_1 \rangle \bar{\nu}_1.$$

Since the action of $Sp(1)$ on $\mathbb{R}^4$ by conjugation is orthogonal and leaves the real axis pointwise fixed (see [10]), the angle $\mathbb{A}_\mathbb{H}(x)$ is independent of the choice of the lifts.

**Proposition 3.2** *Cartan angular invariant $\mathbb{A}_\mathbb{H}(x_1, x_2, x_3)$ is invariant under a permutation of points $x_i$.*

**Proof:** $\langle x_2, x_1, x_3 \rangle = \langle x_2, x_1 \rangle \langle x_1, x_3 \rangle \langle x_3, x_2 \rangle = \overline{\langle x_1, x_2 \rangle} \, \overline{\langle x_3, x_1 \rangle} \, \overline{\langle x_2, x_3 \rangle} = \overline{\langle x_2, x_3 \rangle \langle x_3, x_1 \rangle \langle x_1, x_2 \rangle} = \overline{\langle x_2, x_3, x_1 \rangle}$. The angle between $e_1$ and $\overline{\langle x_2, x_3, x_1 \rangle}$ is the same with the one between $e_1$ and $\langle x_2, x_3, x_1 \rangle$. But also the angle is unchanged under the conjugation by $\langle x_1, x_2 \rangle$, so the angle between $e_1$ and $\langle x_2, x_3, x_1 \rangle$ is equal to the one between $e_1$ and $\langle x_1, x_2 \rangle \langle x_2, x_3, x_1 \rangle \langle x_1, x_2 \rangle^{-1} = \langle x_1, x_2, x_3 \rangle$. So $\mathbb{A}_{\mathbb{H}}(x_1, x_2, x_3) = \mathbb{A}_{\mathbb{H}}(x_2, x_3, x_1) = \mathbb{A}_{\mathbb{H}}(x_2, x_1, x_3)$. ∎

Here we list several properties of the Cartan angular invariant for quaternionic hyperbolic space.

**Proposition 3.3** *Let $x = (x_1, x_2, x_3)$ and $y = (y_1, y_2, y_3)$ be pairs of distinct triples of points in $\overline{H^n_{\mathbb{H}}} = H^n_{\mathbb{H}} \cup \partial H^n_{\mathbb{H}}$. Then $\mathbb{A}_{\mathbb{H}}(x) = \mathbb{A}_{\mathbb{H}}(y)$ if and only if there is an isometry $f \in PSp(n, 1)$ such that $f(x_i) = y_i$ for $i = 1, 2, 3$.*

**Proof:** Applying a homothety (i.e. an element in $PSp(n, 1)$), we may assume that our triples have Hermitian products $X = \langle \tilde{x}_1, \tilde{x}_2, \tilde{x}_3 \rangle$ and $Y = \langle \tilde{y}_1, \tilde{y}_2, \tilde{y}_3 \rangle$ with $|X| = |Y|$.

Now let us assume that $\mathbb{A}_{\mathbb{H}}(x) = \mathbb{A}_{\mathbb{H}}(y)$. This and $|X| = |Y|$ imply that there is an orthogonal transformation $M \in SO(3) \times \{\text{id}\}$ in $\mathbb{H} = \mathbb{R}^4$ that leaves invariant the real axis in $\mathbb{H}$ and maps $X$ to $Y$. Since the conjugation action of $Sp(1)$ in $\mathbb{H}$ is $SO(3)$ action (see §2.1), there is $\mu \in Sp(1)$ such that

$$\langle \tilde{x}_1, \tilde{x}_2, \tilde{x}_3 \rangle = \mu \langle \tilde{y}_1, \tilde{y}_2, \tilde{y}_3 \rangle \bar{\mu} \,.$$

To finish the proof it is enough to choose lifts $\tilde{x}_i$ and $\tilde{y}_i$ of points $x_i$ and $y_i$, $i = 1, 2, 3$, so that $\langle \tilde{x}_i, \tilde{x}_j \rangle = \langle \tilde{y}_i, \tilde{y}_j \rangle$. Indeed, then there is $A \in Sp(n, 1)$ such that $A(\tilde{x}_i) = \tilde{y}_i$, $i = 1, 2, 3$. Then it descends to an element $f \in PSp(n, 1)$ such that $f(x_i) = y_i$ for $i = 1, 2, 3$.

To obtain those lifts, we first replace $\tilde{y}_1$ by $\mu \tilde{y}_1$ (still denote it by $\tilde{y}_1$) and get $\langle \tilde{x}_1, \tilde{x}_2 \rangle \langle \tilde{x}_2, \tilde{x}_3 \rangle \langle \tilde{x}_3, \tilde{x}_1 \rangle = \langle \tilde{y}_1, \tilde{y}_2 \rangle \langle \tilde{y}_2, \tilde{y}_3 \rangle \langle \tilde{y}_3, \tilde{y}_1 \rangle$. Replacing $\tilde{x}_2$ and $\tilde{x}_3$ by $\mu_2 \tilde{x}_2$ and $\mu_3 \tilde{x}_3$ if necessary, we can make $\langle \tilde{x}_2, \tilde{x}_3 \rangle = \langle \tilde{y}_2, \tilde{y}_3 \rangle$ and $\langle \tilde{x}_3, \tilde{x}_1 \rangle = \langle \tilde{y}_3, \tilde{y}_1 \rangle$. Now the equation becomes

$$\langle \tilde{x}_1, \tilde{x}_2 \rangle \langle \tilde{x}_2, \tilde{x}_3 \rangle \langle \tilde{x}_3, \tilde{x}_1 \rangle = |\mu_2|^2 |\mu_3|^2 \langle \tilde{y}_1, \tilde{y}_2 \rangle \langle \tilde{y}_2, \tilde{y}_3 \rangle \langle \tilde{y}_3, \tilde{y}_1 \rangle \,,$$

and we get $\langle \tilde{x}_1, \tilde{x}_2 \rangle = r \langle \tilde{y}_1, \tilde{y}_2 \rangle$ where $r = |\mu_2||\mu_3|$. Then replacing $\tilde{x}_1$, $\tilde{x}_2$, $\tilde{x}_3$ and $\tilde{y}_1$ by $r^{-1} \tilde{x}_1$, $r^{-1} \tilde{x}_2$, $r \tilde{x}_3$ and $r^2 \tilde{y}_1$ respectively, we finally get $\langle \tilde{x}_i, \tilde{x}_j \rangle = \langle \tilde{y}_i, \tilde{y}_j \rangle$, and hence a desired $f \in PSp(n, 1)$.

The converse is trivial. ∎

**Theorem 3.4** *For distinct points $x_1, x_2, x_3 \in \partial H_{\mathbb{H}}^n$, let $\sigma_{12}$ and $\Sigma_{12}$ be real and quaternionic geodesics containing two points $x_1$ and $x_2$, and $\Pi : H_{\mathbb{H}}^n \to \Sigma_{12}$ be the orthogonal projection. Then*

$$|\tan \mathbb{A}_{\mathbb{H}}(x)| = \sinh(d(\Pi x_3, \sigma_{12}))$$

*where $d$ is the hyperbolic distance in $H_{\mathbb{H}}^n$.*

**Proof:** Up to an isometry (in the unit ball model of $H_{\mathbb{H}}^n$), we may assume that the triple $x$ consists of $x_1 = (0, -1)$, $x_2 = (0, 1)$, and $x_3 = (z', z_n)$, whose lifts are $\tilde{x}_1 = (0, -1, 1)$, $\tilde{x}_2 = (0, 1, 1)$, and $\tilde{x}_3 = (z', z_n, 1)$. In this setting $\sigma_{12} = \{(0, t) : t \in \mathbb{R}, |t| < 1\}$, $\Sigma_{12} = \{(0, z) : z \in \mathbb{H}, |z| < 1\}$, and $\langle \tilde{x}_1, \tilde{x}_2, \tilde{x}_3 \rangle = 2(\bar{z}_n - 1)(1 + z_n)$. So we get

$$\bigl|\tan \mathbb{A}_{\mathbb{H}}(x)\bigr| = \frac{|2 \operatorname{Im}(z_n)|}{1 - |z_n|^2}.$$

On the other hand, we note that $\Pi(x_3) = z_n$ and $\Sigma_{12}$ has the Poincaré ball model geometry of $H_{\mathbb{R}}^4$ with sectional curvature $-1$. Choose a hyperbolic two plane in $\Sigma_{12}$ that contains the geodesic $\sigma_{12}$ and $z_n$. This plane is a Poincaré disk with curvature $-1$, where we can write $z_n = \operatorname{Re} z_n + i|\operatorname{Im} z_n|$. Let $d$ be the hyperbolic distance between the point $z_n$ and the real axis in that Poincaré disk. Then a direct calculation shows that $\sinh(d) = |2 \operatorname{Im}(z_n)|/(1 - |z_n|^2)$. ■

**Theorem 3.5** *Three points in the ideal boundary of quaternionic n-space are in the boundary of a real hyperbolic two plane if and only if $\mathbb{A}_{\mathbb{H}}(x) = 0$.*

**Proof:** ($\Rightarrow$) Applying an isometry, we can make these points in the unit ball model of $H_{\mathbb{H}}^n$ as $x_1 = (0, -1)$, $x_2 = (0, 1)$, $x_3 = (0, r)$, where $r$ is real. Then $\langle \tilde{x}_1, \tilde{x}_2, \tilde{x}_3 \rangle$ is real. So the angle $\mathbb{A}_{\mathbb{H}}(x)$ is zero.

($\Leftarrow$) As in the proof of Theorem 3.4, we may assume that $x_1 = (0, -1)$, $x_2 = (0, 1)$, and $x_3 = (z', z_n)$, with lifts $\tilde{x}_1 = (0, -1, 1)$, $\tilde{x}_2 = (0, 1, 1)$, and $\tilde{x}_3 = (z', z_n, 1)$ where $z_n$ must be real in order to have $\mathbb{A}_{\mathbb{H}}(x) = 0$.

Now we use the projection (3) from $N \cup \infty = \partial(\mathbb{H}^{n-1} \times \operatorname{Im} \mathbb{H} \times (0, \infty))$ to the boundary of the unit ball in $\mathbb{H}^n$, where $N$ is the nilpotent group in the Iwasawa decomposition of the automorphism group $PSp(n, 1)$ of $H_{\mathbb{H}}^n$, to see that $x_3$ lies in the boundary of a real hyperbolic two plane.

For real $z_n$, the image $[z, t] \in \mathbb{H}^{n-1} \times \operatorname{Im} \mathbb{H}$ has $t = 0$. So $x_3$ lies on the boundary of a real hyperbolic two plane containing $(0, -1)$ and $(0, 1)$. ■

**Theorem 3.6** *A triple $x = (x_1, x_2, x_3) \in (\partial H_{\mathbb{H}}^n)^3$ lies on the boundary of an $\mathbb{H}$-line if and only if $\mathbb{A}_{\mathbb{H}}(x) = \frac{\pi}{2}$.*

**Proof:** ($\Rightarrow$) If a triple $x = (x_1, x_2, x_3)$ lies on the boundary of an $\mathbb{H}$-line, we can (isometrically) put these three points in the following positions (in the ball model):
$$x_1 = (0, -1), \ x_2 = (0, 1), \ x_3 = (0, z),$$
and lift them to $\tilde{x}_1 = (0, -1, 1)$, $\tilde{x}_2 = (0, 1, 1)$, $\tilde{x}_3 = (0, z, 1)$. Then we have $\langle \tilde{x}_1, \tilde{x}_2, \tilde{x}_3 \rangle = 2(|z|^2 - z + \bar{z} - 1)$. After that by using $|z| = 1$, we get that this number is pure imaginary, and hence $\mathbb{A}_{\mathbb{H}}(x) = \frac{\pi}{2}$.

($\Leftarrow$) Let us assume that $\tilde{x}_1 = (0, -1, 1)$, $\tilde{x}_2 = (0, 1, 1)$ and $\tilde{x}_3 = (z', z_n, 1)$. Then $\langle \tilde{x}_1, \tilde{x}_2, \tilde{x}_3 \rangle = 2(|z_n|^2 - 1 - z_n + \bar{z}_n)$. To get $\langle \tilde{x}_1, \tilde{x}_2, \tilde{x}_3 \rangle$ equal to a pure imaginary number, it is necessary to have $|z_n| = 1$, which implies that $x_3 = (0, z_n)$. This shows that the triple $x$ lies on the $\mathbb{H}$-line $\{(0, \mathbb{H})\}$. ∎

## 4 Toledo invariant in quaternionic space

**Theorem 4.1** *There is a quaternionic Kähler four form $\omega$ such that its restriction to any $\mathbb{H}$-line is its volume form. Furthermore it is a closed form and its evaluation on four vectors two of which span a totally real geodesic two plane is zero.*

**Proof:** Since we may identify the quaternionic hyperbolic $n$-space with $PSp(n,1)/Sp(n)Sp(1)$, it suffices to prove the existence of a 4-form that is invariant under $Sp(n)Sp(1)$ at the tangent space of the quaternionic space at the origin, which is $\mathbb{R}^{4n}$. The action of $Sp(n)Sp(1)$ on the tangent space is given by $v \to Sp(1)^{-1}vSp(n)$. Such a 4-form $\omega$ is given in [12] (Proposition 1.1). By the construction it is invariant under the isometry group and so is a parallel 4-form. Then due to the standard fact in Riemannian geometry, this form is closed. To show the last statement of the theorem, recall the construction of the form $\omega$. Let $L$ be a quaternionic line and $\nu_L$ a volume form on $L$. Let $\Pi_L$ be the orthogonal projection onto $L$. Then
$$\Omega = \int_{\mathbb{H}P^{n-1}} \Pi_L^* \nu_L dL$$
is an average over $Sp(n)$, and finally $\omega$ is its transport by $PSp(n,1)$. Finally to show that evaluation of $\omega$ on four vectors two of which span a totally real geodesic space $H$ is identically zero, one should note the following. Let $H$ be $\{(\mathbb{R}, \mathbb{R}, 0, \cdots, 0)\}$ a standard totally real 2-space in the unit ball model and $e_1, e_2$ be its standard basis. If $L$ is a quaternionic line, let $A \in Sp(n)$ be such that $LA = \{(\mathbb{H}, 0, \cdots, 0)\}$ is the first coordinate line $L_1$. Then

$$\Pi_L{}^*\nu_L(e_1, e_2, v, w) = \nu_{L_1}(\Pi_{L_1}(e_1 A), \Pi_{L_1}(e_2 A), \Pi_{L_1}(vA), \Pi_{L_1}(wA)).$$

Let $A'$ be an element in $Sp(n)$ obtained from $A$ by interchanging the first and the second rows. Set $L' = L_1 A'^{-1}$. Then

$$\Pi_{L'}{}^*\nu_{L'}(e_1, e_2, v, w) = \nu_{L_1}(\Pi_{L_1}(e_1 A'), \cdots, \Pi_{L_1}(wA'))$$

$$= \nu_{L_1}(\Pi_{L_1}(e_2 A), \Pi_{L_1}(e_1 A), \Pi_{L_1}(vA), \Pi_{L_1}(wA)).$$

So $\Pi_L{}^*\nu_L(e_1, e_2, v, w) + \Pi_{L'}{}^*\nu_{L'}(e_1, e_2, v, w) = 0$. This way the average over the $\mathbb{H}$-lines through the origin is zero, and $\omega$ satisfies the desired property. ∎

Let $x = (x_1, x_2, x_3)$ be a triple of distinct points in $H_\mathbb{H}^n \cup \partial H_\mathbb{H}^n$, and let $\sigma_{12}$ and $\Sigma_{12}$ be unique real and quaternionic geodesics containing $x_1$ and $x_2$. Then the stabilizer of the geodesic $\sigma_{12}$ is $Sp(1) \times Sp(n-1) \times \mathbb{R}$. Let $\Pi \colon H_\mathbb{H}^n \to \Sigma_{12}$ be the orthogonal projection onto the $\mathbb{H}$-line $\Sigma_{12}$ and $\triangle \subset \Sigma_{12}$ be a triangle with geodesic edges spanned by vertices $\Pi(x_1) = x_1$, $\Pi(x_2) = x_2$ and $\Pi(x_3)$. We define $\triangle_{\Pi(x)}$ as the $Sp(1)$-orbit of the triangle $\triangle$. This is a four dimensional object in the $\mathbb{H}$-line. Finally we define $\triangle_x \subset H_\mathbb{H}^n$ to be any four dimensional object containing our given triple $x = (x_1, x_2, x_3)$ and whose $\Pi$-projection on the $\mathbb{H}$-line is $\triangle_{\Pi(x)}$.

**Definition 4.2** *The quaternionic Toledo invariant $\tau_\mathbb{H}(x)$ of a triple $x$ is*

$$\tau_\mathbb{H}(x) = \frac{1}{4\pi} \int_{\triangle_x} \omega$$

*where $\omega$ is a quaternionic Kähler form in Theorem 4.1.*

We shall prove that our definition does not depend on a choice of lift $\triangle_x$.

**Theorem 4.3** *For any triple $x = (x_1, x_2, x_3) \in (\overline{H_\mathbb{H}^n})^3$, $\tau_\mathbb{H}(x) = 2\mathbb{A}_\mathbb{H}(x)$.*

**Proof:** Let $\triangle$ be a 5-dimensional object in $H_{\mathbb{H}}^n$ whose faces consist of $\triangle_x$, $\triangle_{\Pi(x)}$ and $\sigma$ where $\sigma$ is a (totally real geodesic) vertical 4-dimensional surface connecting $\partial \triangle_x$ and $\partial \triangle_{\Pi(x)}$. Due to the Stoke's theorem and $d\omega = 0$, $\int_\triangle d\omega = \int_{\partial \triangle} \omega = 0$. On the other hand, due to Lemma 2.2 the tangent space at any point $z \in \sigma$ contains two vectors spanning a totally geodesic real plane. This and Theorem 4.1 imply that $\int_\sigma \omega = 0$, and we get

$$\int_{\triangle_x} \omega = \int_{\triangle_{\Pi(x)}} \omega = \int_{\triangle_{\Pi(x)}} \omega_{vol},$$

where $\Pi \colon H_{\mathbb{H}}^n \to \Sigma_{12}$ is the orthogonal projection, $\Sigma_{12}$ is an $\mathbb{H}$-line containing $x_1$ and $x_2$, and $\omega_{vol}$ is a volume form on $\Sigma_{12}$. Let $x_1 = (0, ..., 0, -1)$, $x_2 = (0, ..., 0, 1)$, and $x_3 = (z_1, \cdots, z_n)$. Then $\int_{\Pi(\triangle_x)} \omega_{vol}$ is the integral of the hyperbolic volume form on $H_{\mathbb{R}}^4$ over $\Pi(\triangle_x)$. Here $\Pi(\triangle_x)$ is the $Sp(1)$-orbit of the geodesic triangle $\triangle = \triangle(x_1, x_2, z_n)$ in $H_{\mathbb{R}}^4$ with vertices $x_1$, $x_2$, and $\Pi(x_3) = (0, ..., 0, z_n)$. But the hyperbolic metric on $H_{\mathbb{R}}^4$ along the geodesic $\sigma_{12}$ can be written as

$$d\rho^2 + \sinh^2 \rho \cdot d\theta^2 + \cosh^2 \rho \cdot d\lambda^2,$$

where $\rho$ is a hyperbolic distance from $\sigma_{12}$, $\lambda$ is a distance along $\sigma_{12}$ and $\theta$ is a rotation parameter around $\sigma_{12}$. Then

$$\int_{\triangle_{\Pi(x)}} \omega_{vol} = \int_{Sp(1)\triangle(x_1,x_2,z_n)} \omega_{vol} = \int_{S^2} \int_{\triangle(x_1,x_2,z_n)} \sinh \rho \cosh \rho \, d\rho \, d\lambda \, d\theta$$

$$= \text{Area}(\triangle(x_1, x_2, z_n)) \, \text{Area}(S^2) = 4\pi \, \text{Area}(\triangle(x_1, x_2, z_n)).$$

Now it suffices to show that $\text{Area}(\triangle(x_1, x_2, z_n)) = 2\mathbb{A}(x)$. Indeed, the Cartan angular invariant $\mathbb{A}_{\mathbb{H}}(x)$ is the angle between the line $\mathbb{R}e_1 \subset \mathbb{R}^4$ and the vector representing a quaternion $(\bar{z}_n - 1)(1 + z_n)$ because the Hermitian triple product $\langle \tilde{x}_1, \tilde{x}_2, \tilde{x}_3 \rangle = 2(\bar{z}_n - 1)(1 + z_n)$.

By choosing a hyperbolic two plane containing $z_n$, we may assume that the points $x_1, x_2$ and $z_n$ lie on the Poincaré unit disk, and $z_n = \text{Re } z_n + i |\text{Im } z_n|$ is a complex number. Then $\mathbb{A}_{\mathbb{H}}(x)$ becomes $|\arg((1 - \bar{z}_n)(1 + z_n))|$. Using the Cayley transformation $\phi$ from the Poincaré disk to the right half plane, we have $\phi(z_n) = \frac{1}{2}(1 - z_n)/(1 + z_n)$. It is easy to see that in this model of the hyperbolic plane, the area of the triangle $\triangle_{0,\infty,\phi(z_n)}$ is $2|\arg \phi(z_n)|$. But $|\arg \phi(z_n)| = |\arg((1 - z_n)(1 + \bar{z}_n))| = \mathbb{A}_{\mathbb{H}}(x)$. This finishes the proof. ∎

Since the angular invariant $\mathbb{A}_{\mathbb{H}}(x)$ is invariant under the permutation of the variables due to Proposition 3.2, the above theorem immediately implies:

**Corollary 4.4** *The quaternionic Toledo invariant $\tau_{\mathbb{H}}(x)$ of a triple $x$, $x = (x_1, x_2, x_3) \in (\overline{H^n_{\mathbb{H}}})^3$, is invariant under permutations of the variables.*

In other words, the above fact shows that the integral of the quaternionic Kähler form in the definition of Toledo invariant is independent of the choice of a geodesic edge $\overline{x_i x_j}$ used in constructing $\triangle_x$.

## 5 Angular invariant in the Cayley plane

In the Cayley hyperbolic plane, there is a closed 8-form $\omega$ satisfying the same properties as the quaternionic 4-form in Theorem 4.1, cf. [11]. Taking into account Theorem 4.3, we can use this closed 8-form $\omega$ for the definition of an angular invariant in the Cayley plane. Since the subgroup of $F_4^{-20}$ fixing $\infty$ and the origin of the Carnot group $N \cong \partial H^2_{\mathbb{O}} \setminus \{\infty\}$ (in the Iwasawa decomposition $F_4^{-20} = KAN$) is $Spin(7)A$, the stabilizer of a real geodesic in $H^2_{\mathbb{O}}$ is conjugate to $Spin(7)$. For a triple $x = (x_1, x_2, x_3)$ of points on the ideal boundary $\partial H^2_{\mathbb{O}}$, let $\Pi \colon H^2_{\mathbb{O}} \to \Sigma_{12}$ be the orthogonal projection to the unique octonionic geodesic ($\mathbb{O}$-line) $\Sigma_{12}$ spanned by $x_1$ and $x_2$. We define $\triangle_x$ to be the $Spin(7)$-orbit $Spin(7)\triangle$ of $\triangle$. Here $Spin(7)$ is the stabilizer of the real geodesic connecting the points $x_1$ and $x_2$, and $\triangle$ is the geodesic 2-dimensional triangle (in the $\mathbb{O}$-line $\Sigma_{12}$) with the vertices $x_1$, $x_2$, and $\Pi(x_3)$. Then $\triangle_x$ is an 8-dimensional object, and we define:

**Definition 5.1** *The octonionic angular invariant $\mathbb{A}_{\mathbb{O}}(x)$ of a triple $x$ is*

$$\mathbb{A}_{\mathbb{O}}(x) = \frac{1}{\text{Vol}(S^6)} \int_{\triangle_x} \omega$$

*where $\omega$ is the octonionic Kähler form.*

From this definition and properties of the orthogonal projection $\Pi$ of $H^2_{\mathbb{O}}$ onto an $\mathbb{O}$-line we immediately have:

**Theorem 5.2** *Let $x = (x_1, x_2, x_3)$ be a triple of points on the ideal boundary $\partial H^2_{\mathbb{O}}$. Then its angular invariant $\mathbb{A}_{\mathbb{O}}(x)$ has the following properties:*

1. $\mathbb{A}_\mathbb{O}$ is invariant under isometries of $H^2_\mathbb{O}$, so triples $x$ and $x'$ with different values $\mathbb{A}_\mathbb{O}(x) \neq \mathbb{A}_\mathbb{O}(x')$ are not isometrically equivalent;

2. $\mathbb{A}_\mathbb{O}(x) = 0$ if and only if a triple $x$ lies on the boundary of a totally real hyperbolic 2-plane in $H^2_\mathbb{O}$;

3. $\mathbb{A}_\mathbb{O}(x) = \pi/2$ if and only if $x$ lies on the boundary of an $\mathbb{O}$-line.

As for quaternionic case, we note that the octonionic angular invariant $\mathbb{A}_\mathbb{O}(x)$ cannot distinguish different shapes of $k$-spheres inside the sphere at infinity of an $\mathbb{O}$-line, in particular, a $k$-dimensional linear subspace $L$ and the result $L'$ of its bending along a hyperplane, see next Section.

## 6 Bendings in rank one symmetric spaces

Let $M_0$ be a locally symmetric space of rank one with the fundamental group (in the orbifold category) $\Gamma = \pi_1(M_0)$ which acts discretely in its universal cover $\widetilde{M_0}$, i.e. in a hyperbolic $\mathbb{F}$-space $H^n_\mathbb{F}$. The deformation (Teichmüller) space $\mathcal{T}(M_0)$ can be defined as the space of isotopy classes of marked $\mathbb{F}$-hyperbolic structures on $M_0$. By using holonomy representations induced by the holonomy maps of structures on $M_0$, one can study instead the variety $\mathcal{T}(\Gamma)$ of conjugacy classes of discrete faithful representations $\rho\colon \Gamma \to \mathrm{Isom}\, H^n_\mathbb{F}$. Here $\mathcal{T}(\Gamma) = \mathcal{R}_0(\Gamma)/\mathrm{Isom}\, H^n_\mathbb{F}$, and the variety $\mathcal{R}_0(\Gamma) \subset \mathrm{Hom}(\Gamma, \mathrm{Isom}\, H^n_\mathbb{F})$ consists of discrete faithful representations $\rho$ of the group $\Gamma$ whose co-volume, $\mathrm{Vol}(H^n_\mathbb{F}/\Gamma)$, may be infinite.

Here we shall present a canonical construction of (bending) deformations of a wide class of analytic $\mathbb{F}$-hyperbolic manifolds $M_0$. This will show non-triviality of the Teichmüller space $\mathcal{T}(M_0)$ for such manifolds, despite many rigidity results [23, 22, 12]. First such deformations were given by the first author for real hyperbolic n-manifolds, $n > 3$, and after Thurston's "Mickey Mouse" example, were called bendings of $M_0$ along its totally geodesic hypersurfaces, see [1]. Recently, the first author [3] defined bending deformations for complex hyperbolic manifolds where such bendings are associated to (real) closed geodesics in $M_0$ and can be induced by equivariant quasiconformal homeomorphisms.

However, quaternionic and octonionic hyperbolic structures differ from such real and complex ones. Indeed they are much more rigid due to Corlette's [12] and Gromov-Schoen's [19] superrigidity of their uniform lattices.

## 6.1 Bending representations

Here we shall show to what extent bending deformations really work for quaternionic and octonionic hyperbolic structures on $M_0 = H_\mathbb{F}^n/\Gamma$, or when do they provide paths in the Teichmüller space $\mathcal{T}(\Gamma)$ that consist of conjugacy classes of non-trivial quasi-Fuchsian representations $\rho_t \colon \Gamma \to \text{Isom}\, H_\mathbb{F}^n$.

We start with pointing out some situations when such bendings are impossible. First, let $\Gamma \subset \text{Isom}\, H_\mathbb{F}^n$ be a discrete group preserving a quaternionic or Cayley line $H_\mathbb{F}^1 \subset H_\mathbb{F}^n$. The induced metric on $H_\mathbb{F}^1$ is the Poincaré metric of $H_\mathbb{R}^k$ where $k$ equals either four or eight respectively. So if the group $\Gamma$ has a bending hyperplane $P \subset H_\mathbb{F}^1$, $P \cong H_\mathbb{R}^{k-1}$, one can bend $\Gamma$ inside a bigger *real* hyperbolic space $H_\mathbb{R}^{k+1}$, see [1]. Nevertheless, each analytic rotation $U$ in $H_\mathbb{F}^n$ pointwise fixing $P$ would be the identity on the whole $\mathbb{F}$-line, which makes impossible a naive attempt to bend the group $\Gamma$ inside the whole space $H_\mathbb{F}^n$, cf. (10). Second "non-bending" situation is for lattices $\Gamma \subset \text{Isom}\, H_\mathbb{F}^k$, $2 \leq k \leq n$. Here we see that the set of fixed points of a compact one parameter Lie subgroup in $\text{Isom}\, H_\mathbb{F}^n$ (whose elements one should use for bending as in (10)) is an analytic geodesic subspace, so its real codimension is at least two. So we cannot use such non-separating planes as bending ones.

Let $\Gamma \subset \text{Isom}\, H_\mathbb{R}^n \subset \text{Isom}\, H_\mathbb{F}^n$ be a discrete group (not necessarily a lattice) acting by isometries in a totally geodesic real subspace $H_\mathbb{R}^n \subset H_\mathbb{F}^n$ as well as in the whole $\mathbb{F}$-space $H_\mathbb{F}^n$ with a distance function $d$. The main condition on the group $\Gamma$ is that it has a *real bending hyperplane* $P \cong H_\mathbb{R}^{n-1} \subset H_\mathbb{R}^n$:

$$\inf\{d(\gamma(P), P) \colon \gamma \in \Gamma \backslash \Gamma_P\} = d_P > 0, \qquad (8)$$

where we additionally assume that the stabilizer $\Gamma_P = \{\gamma \in \Gamma : \gamma(P) = P\}$ of $P$ in $\Gamma$ is a proper subgroup, $\Gamma_P \neq \Gamma$.

Due to the tubular neighborhood theorem [9, 8], any real hyperbolic $n$-manifold with an embedded closed totally geodesic hypersurface $S$ of area $A$ has a tubular neighborhood of $S$ of width $c_n(A)$ that depends only on dimension $n$ and area $A$. This implies that the above *real bending hyperplane* condition on the group $\Gamma \subset \text{Isom}\, H_\mathbb{R}^n$ with a uniform lattice $\Gamma_P \subset \text{Isom}\, H_\mathbb{R}^{n-1}$ is equivalent to the existence of an embedded totally geodesic hyperspace $P/\Gamma_P \subset H_\mathbb{R}^n/\Gamma$, which is called *bending suborbifold*. In other words, we shall be able to make our bending deformations of a $\mathbb{F}$-hyperbolic manifold $M_0$ if it is homotopy equivalent to a totally real submanifold with embedded totally geodesic hypersurface. In the unit ball models, we have:

$$P = H_\mathbb{R}^{n-1} \cong B_\mathbb{R}^{n-1} \subset B_\mathbb{R}^n \cong H_\mathbb{R}^n \subset H_\mathbb{F}^n \cong B_\mathbb{F}^n.$$

Also we have a $\mathbb{F}$-hyperplane $H_\mathbb{F}^{n-1} \subset H_\mathbb{F}^n$ that orthogonally intersects our real subspace $H_\mathbb{R}^n$ along $P$. At a point $x_0 \in P$, the orthogonal complement to $P$ in $H_\mathbb{R}^n$ is a real geodesic $\ell \subset H_\mathbb{R}^n$. This geodesic $\ell$ is the (orthogonal) intersection of $H_\mathbb{R}^n$ and the $\mathbb{F}$-geodesic $H_\mathbb{F}^1 \subset H_\mathbb{F}^n$ that is the orthogonal complement to the $\mathbb{F}$-hyperplane $H_\mathbb{F}^{n-1}$ at $x_0$. The $\mathbb{F}$-hyperplane $H_\mathbb{F}^{n-1}$ is the fixed point set of a connected compact Lie subgroup $\mathcal{U}_\mathbb{F}$ of automorphisms of $H_\mathbb{F}^n$. It corresponds to the group of unit $\mathbb{F}$-numbers and acts in $H_\mathbb{F}^n$ and $H_\mathbb{F}^1$ by orthogonal rotations. Its real dimension $k_\mathbb{F} = \dim \mathcal{U}_\mathbb{F}$ is 1, 3 and 7 respectively for complex, quaternionic and octonionic cases. Exept involutions, all other elements $U \in \mathcal{U}$ map the geodesic $\ell$ to geodesics $U(\ell)$ that entersect $\ell$ at $x_0$ and do not lie in $H_\mathbb{R}^n$.

In horospherical coordinates, $\mathbb{F}^{n-1} \times \operatorname{Im} \mathbb{F} \times [0, \infty)$, by looking at infinity of all those spaces we see that any our rotation $U \in \mathcal{U}$ pointwise fixes a subspace $\overline{\mathbb{R}^{n-2}} = \partial H_\mathbb{R}^{n-1} \subset \partial H_\mathbb{F}^{n-1}$ that lies in the real subspace $\overline{\mathbb{R}^{n-1}} = \partial H_\mathbb{R}^n$ of the *horizontal* plane $\mathbb{F}^{n-1} \times \{0\}$ in the Carnot group $\mathbb{F}^{n-1} \times \operatorname{Im} \mathbb{F}$.

Now let $\Gamma \subset PO(n,1) \subset \operatorname{Isom} H_\mathbb{F}^n$ be a discrete group preserving a totally real n-space $H_\mathbb{R}^n \subset H_\mathbb{F}^n$ and satisfying (8), that is having a real bending hyperplane $P \subset H_\mathbb{R}^n$ with a stabilizer $\Gamma_P \subset \Gamma$, $\Gamma_P \neq \Gamma$. Then, depending on whether the suborbifold $P/\Gamma_P$ is separating or not in $H_\mathbb{R}^n/\Gamma$, the group $\Gamma$ is either a free amalgamated product or $HNN$-group:

$$\Gamma = \Gamma_1 *_{\Gamma_P} \Gamma_2 \quad \text{or} \quad \Gamma = \Gamma_1 *_{\Gamma_P} = \langle \Gamma_1, \gamma_2 \rangle, \qquad (9)$$

where $\Gamma_1$ and $\Gamma_2$ are subgroups of $\Gamma$, $\Gamma_1 \cap \Gamma_2 = \Gamma_P$, and $\gamma_2 \in \Gamma \backslash \Gamma_1$ conjugates two subgroups of $\Gamma_1$ that are isomorphic to $\Gamma_P$. We shall describe our bending deformation of $\Gamma$ in $\operatorname{Isom} H_\mathbb{F}^n$, along the real bending hyperplane $P \subset H_\mathbb{R}^n$.

To do that, we consider a one parameter Lie subgroup $\{U_\eta\}$ of the above compact Lie group $\mathcal{U} \subset \operatorname{Isom} H_\mathbb{F}^n$ that pointwise fixes the $\mathbb{F}$-hyperplane $H_\mathbb{F}^{n-1}$. Then, in both cases of free amalgamated products or $HNN$-groups in (9), we define deformed groups $\Gamma_\eta \subset \operatorname{Isom} H_\mathbb{F}^n$ correspondingly as

$$\Gamma_\eta = \Gamma_1 *_{\Gamma_P} U_\eta \Gamma_2 U_\eta^{-1} \quad \text{or} \quad \Gamma_\eta = \langle \Gamma_1, U_\eta \gamma_2 \rangle \cong \Gamma_1 *_{\Gamma_P}. \qquad (10)$$

We shall observe later that, for small enough parameters $\eta$, the deformed groups $\Gamma_\eta \subset \operatorname{Isom} H_\mathbb{F}^n$ are discrete. In general situation, this discreteness follows from a topological conjugation of $\Gamma_\eta$ and $\Gamma$ in the whole space $H_\mathbb{F}^n$. For convex co-compact torsion free groups $\Gamma_\eta \subset H_\mathbb{F}^n$ one can also use a result by Chen and Greenberg (see Cor. 4.5.5 in [10]) that such a pure loxodromic group is indeed a discrete group.

Furthermore, a direct application of Cartan angular (or Toledo) invariants implies that, for non-trivial subgroups $\Gamma_i$ (which is the case for real lattices $\Gamma$), the deformed groups $\Gamma_\eta$ are non-trivial quasi-Fuchsian groups in $H_\mathbb{F}^n$. The limit set of $\Gamma_\eta$ is then a Jordan $(n-1)$-dimensional sphere in $\partial H_\mathbb{F}^n$ which has the Hausdorff dimension strictly greater than $n-1$.

Also we see from (10) that the family of bending representations

$$\rho_\eta \colon \Gamma \to \Gamma_\eta \subset \operatorname{Isom} H_\mathbb{F}^n$$

analitically depends on the parameter $\eta$, so it defines smooth $k_\mathbb{F}$-dimensional curves in Teichmüller spaces $\mathcal{T}(\Gamma)$ and $\mathcal{T}(M_0)$, where $k_\mathbb{F} = \dim \mathcal{U}_\mathbb{F}$ is 1, 3 and 7 respectively for complex, quaternionic and octonionic cases.

## 6.2 Bending of a surface group in $\operatorname{Isom} H_\mathbb{F}^n$

Now we start with a totally real geodesic surface $S = H_\mathbb{R}^2/\Gamma$ in the $\mathbb{F}$-hyperbolic manifold $M_0 = H_\mathbb{F}^2/\Gamma$, where $\Gamma \subset PO(2,1) \subset \operatorname{Isom} H_\mathbb{F}^2$ is a given discrete group, and fix a simple geodesic $\alpha$ on $S$. Here we may obviously change $H_\mathbb{F}^2$ to $n$-dimensional $\mathbb{F}$-hyperbolic space, in particular to quaternionic space $H_\mathbb{H}^n$, which provides more parameters for our bending deformation because of more freedom of rotations around $\mathbb{F}$-line in $n$-space. Let $P \subset H_\mathbb{R}^2 \subset H_\mathbb{F}^2$ be a real geodesic covering the loop $\alpha$ and with ends (in horospherical coordinates) at $\infty$ and the origin of the Carnot group $N = \mathbb{F} \times \operatorname{Im} \mathbb{F}$. Deforming the surface $S$ and its holonomy group $\Gamma$ in $PO(2,1)$, we may assume that the hyperbolic length $\varepsilon$ of $\alpha$ is sufficiently small and the radius $\delta$ of its tubular neighborhood is big enough. It is due to the following claim [3], whose conditions in particular hold for $\varepsilon < 0.791411$ and $\ln 63 \leq \delta < \operatorname{arccsch}\left(2 \sinh \frac{\varepsilon}{4}\right)$:

**Proposition 6.1** *Let $\gamma_\alpha$ be a hyperbolic element of a non-elementary discrete group $\Gamma \subset PO(2,1) \subset \operatorname{Isom} H_\mathbb{F}^2$ with translation length $\varepsilon$ along its axis $P \subset H_\mathbb{R}^2$. Then any tubular neighborhood $U_\delta(P)$ of the axis $P$ of radius $\delta > 0$ is precisely invariant with respect to its stabilizer $\Gamma_P \subset \Gamma$ if*

$$\sinh\left(\frac{\varepsilon}{4}\right) \cdot \sinh\left(\frac{\delta}{2}\right) \leq \frac{1}{2}. \tag{11}$$

*Furthermore, for sufficiently small $\varepsilon$, $\varepsilon < 4\delta$, the Dirichlet polyhedron $D_z(\Gamma)$ of the group $\Gamma$ centered at a point $z \in P$ has two sides $\sigma$ and $\sigma'$ intersecting the axis $P$ and such that $\gamma_\alpha(\sigma) = \sigma'$.*

Now applying Proposition 6.1, the description of bisectors in $H_\mathbb{F}^n$ given in Section 2.2 and Theorem 2.4, we may consider the length $\varepsilon$ of the geodesic $\alpha$ so small that closures of all equidistant halfspaces in $\overline{H_\mathbb{F}^2}\setminus\{\infty\}$ bounded by those bisectors and disjoint from the Dirichlet polyhedron $D_z(\Gamma)$ do not intersect the co-vertical bisector whose infinity is $\operatorname{Im}\mathbb{F} \times \operatorname{Im}\mathbb{F} \subset \mathbb{F} \times \operatorname{Im}\mathbb{F}$. Then the group $\Gamma$ and its subgroups $\Gamma_P, \Gamma_1, \Gamma_2$ in (9) have Dirichlet polyhedra centered at a point $z \in P = (0, \infty)$, whose intersections with the hyperbolic 2-plane $H_\mathbb{R}^2$ have the standard shapes for such products (9), see [3].

In particular we have that all bisectors (4) bounding those Dirichlet polyhedra, except two bisectors $\mathfrak{S}$ and $\mathfrak{S}'$ that are identified by the hyperbolic generator $\gamma_\alpha$ of the stabilizer $\Gamma_P$ of the axis $P$, lie in sufficiently small "cone neighborhoods" $C_+$ and $C_-$ of the rays $\mathbb{R}_+$ and $\mathbb{R}_-$ of the real line $\mathbb{R} \times \{0\} = \partial H_\mathbb{R}^2\setminus\{\infty\} \subset \mathbb{F} \times \operatorname{Im}\mathbb{F} = N$. Let us assume that the radii (in the Cygan metric (2)) of the spheres at infinity of the bisectors $\mathfrak{S}$ and $\mathfrak{S}'$ are 1 and $r_0 > 1$, and take a sufficiently small $\varepsilon$, $0 < \varepsilon << r_0 - 1$. Then the cone neighborhoods $C_+, C_- \subset \overline{H_\mathbb{F}^n}\setminus\{\infty\} = \mathbb{F}^{n-1} \times \operatorname{Im}\mathbb{F} \times [0, +\infty)$ are correspondingly the cones of the $\varepsilon$-neighborhoods of the points $(1,0,0), (-1,0,0) \in \mathbb{F}^{n-1} \times \operatorname{Im}\mathbb{F} \times [0, +\infty)$ with respect to the Cygan metric $\rho_c$ in $\overline{H_\mathbb{F}^n}\setminus\{\infty\}$, from the origin.

Now we specify numbers $\eta$ and $\zeta$, $0 < \zeta < \pi/2$, $0 \leq \eta < \pi - 2\zeta$, and such that the intersection $C_+ \cap (\mathbb{F} \times \{0\})$ is contained in the $\zeta$-angle neighborhood $\mathbb{F}_\zeta(\mathbb{R}_+) \subset \mathbb{F}$ of the positive real ray $\mathbb{R}_+$ in $\mathbb{F}$ which is considered as Euclidean space of dimension $\dim \mathbb{F}$. Then for a one parameter Lie subgroup $\{U_t\}$ of the compact Lie group $\mathcal{U} \subset \operatorname{Isom}\mathbb{F}$, we define a bending homeomorphism $\phi = \phi_{\eta,\zeta}\colon \mathbb{F} \to \mathbb{F}$ which bends the real axis $\mathbb{R} \subset \mathbb{F}$ at the origin by the angle $\eta$ and is $\Gamma_P$-equivariant. This homeomorphism $\phi_{\eta,\zeta}$ is defined as follows:

1. It is the identity in the neighborhood $\mathbb{F}_\zeta(\mathbb{R}_-) \subset \mathbb{F}$ of the real ray $\mathbb{R}_-$;

2. It is the rotation $U_\eta$ in the neighborhood $\mathbb{F}_\zeta(\mathbb{R}_+) \subset \mathbb{F}$ of the ray $\mathbb{R}_+$;

3. It preserves each sphere centered at $0 \in \mathbb{F}$ (in Cygan metric);

4. It maps each ray in $\mathbb{F}$ amanating from $0 \in \mathbb{F}$ to a similar ray.

Foliating the punctured Carnot group $N\setminus\{0\}$ and the whole $\mathbb{F}$-hyperbolic space correspondingly by C-spheres $S(0, r) = \{z \in N : |z|_c = r\}$ of radii $r > 0$ and by bisectors having those C-spheres at their infinity, we can extend this bending homeomorphism to a $\Gamma_P$-equivariant homeomorphism of the whole space $\overline{H_\mathbb{F}^2}$ (compatible with its structure, compare our construction in [3]).

However, due to P.Pansu's [24] rigidity, we cannot make them quasiisometric, so we may restrict ourselves to a simpler topological extension. Namely, we extend the cone neighborhoods $\mathbb{F}_\zeta(\mathbb{R}_\pm) \subset \mathbb{F}$ of the real rays $\mathbb{R}_\pm$ to their cone neighborhoods $N_\zeta(\mathbb{R}_\pm) \subset N$ and $H_\zeta(\mathbb{R}_\pm) \subset \overline{H_\mathbb{F}^2}$. These new cones intersect $\mathbb{F}$ (and correspondingly $N$) along the previously defined cones and contain all bisectors (and their spinal spheres at infinity) that bound the Dirichlet polyhedron $D_z(\Gamma)$ and differ from two bisectors $\mathfrak{S}$ and $\mathfrak{S}'$ identified by the generator $\gamma_\alpha \in \Gamma_P$. Furthermore, due to our construction, the closures of these neighborhoods intersect the closure of the co-vertical bisector whose infinity is $\operatorname{Im}\mathbb{F} \times \operatorname{Im}\mathbb{F} \subset \mathbb{F} \times \operatorname{Im}\mathbb{F}$ only at the origin of $N = \mathbb{F} \times \operatorname{Im}\mathbb{F}$.

Now we define an elementary bending homeomorphism $\varphi_{\eta,\zeta} : \overline{H_\mathbb{F}^2} \to \overline{H_\mathbb{F}^2}$ as the following $\Gamma_P$-equivariant extension of the homeomorphism $\phi_{\eta,\zeta}$:

1. It is identity in the closure of the neighborhood $H_\zeta(\mathbb{R}_-) \subset \overline{H_\mathbb{F}^2}$ of the real ray $\mathbb{R}_-$ and in the closure of the $\mathbb{F}$-geodesic $\{0\} \times \operatorname{Im}\mathbb{F} \times (0,\infty) \subset H_\mathbb{F}^2$;

2. It is the rotation $U_\eta$ in the neighborhood $H_\zeta(\mathbb{R}_+) \subset \overline{H_\mathbb{F}^2}$ of the ray $\mathbb{R}_+$;

3. It preserves bisectors $\mathfrak{S}_r$ containing at infinity $\mathbb{F}$-spheres $S(0,r) \subset N$;

4. It commutes with the generator $\gamma_\alpha$ of $\Gamma_P$ on the bisector complements $\overline{\mathfrak{S}}\backslash H_\zeta(\mathbb{R}_\pm)$ and $\overline{\mathfrak{S}'}\backslash H_\zeta(\mathbb{R}_\pm)$.

Such an extension can be made as usual by extending a homeomorphism to a topological ball from its values on the boundary sphere. Meanwhile we see some "isometric behavior" of our homeomorphism $\varphi_{\eta,\zeta}$ on faces of the Dirichlet polyhedron $D_z(\Gamma)$ (except of bisectors $\mathfrak{S}$ and $\mathfrak{S}'$). Namely, let $\mathfrak{S}^+ \subset H_\mathbb{F}^2$ be a "half-space" disjoint from $D_z(\Gamma)$ and bounded by a bisector different from bisectors $\mathfrak{S}_r, r > 0$, and contain a side of the polyhedron $D_z(\Gamma)$. Then there is an open neighborhood $N(\overline{\mathfrak{S}^+}) \subset \overline{H_\mathbb{F}^2}$ such that the restriction of $\varphi_{\eta,\zeta}$ to it is either the identity or coincides with the rotation $U_\eta \subset \operatorname{Isom} H_\mathbb{F}^2$ by the angle $\eta$ about the "vertical" $\mathbb{F}$-geodesic (whose infinity contains the center $\{0\} \times \operatorname{Im}\mathbb{F}$ of $N$).

Clearly, the image $D_\eta = \varphi_{\eta,\zeta}(D_z(\Gamma))$ is a polyhedron in $H_\mathbb{F}^2$ bounded by bisectors. Furthermore, there is a natural identification of its sides induced by $\varphi_{\eta,\zeta}$. Actually the pairs of sides preserved by $\varphi_{\eta,\zeta}$ are identified by the original generators of the group $\Gamma_1 \subset \Gamma$. For other sides $\mathfrak{s}_\eta$ of $D_\eta$, which are images of corresponding sides $\mathfrak{s} \subset D_z(\Gamma)$ under the unitary rotation $U_\eta$, we define side pairings by using the group $\Gamma$ decomposition. For $\Gamma = \Gamma_1 *_{\Gamma_P} \Gamma_2$, we change the original side pairings $\gamma \in \Gamma_2$ of $D_z(\Gamma)$-sides to the $\mathbb{F}$-hyperbolic isometries $U_\eta \gamma U_\eta^{-1}$. For HNN-extension $\Gamma = \Gamma_1 *_{\Gamma_P} = \langle \Gamma_1, \gamma_2 \rangle$, we change the original side pairing $\gamma_2 \in \Gamma$ of $D_z(\Gamma)$-sides to the $\mathbb{F}$-hyperbolic isometry $U_\eta \gamma_2$. So the deformed groups $\Gamma_\eta \subset \operatorname{Isom} H^2_\mathbb{F}$ have the form (10).

This shows that the representations $\Gamma \to \Gamma_\eta \subset \operatorname{Isom} H^2_\mathbb{F}$ do not depend on $\zeta$ and analytically depend on the angle parameter $\eta$. For small enough angle $\eta$, it also follows that the behavior of neighboring polyhedra $\gamma'(D_\eta)$, $\gamma' \in \Gamma_\eta$ bounded by bisectors, around edges of the polyhedron $D_\eta$, is the same as of those $\gamma(D_z(\Gamma))$, $\gamma \in \Gamma$, around edges of the Dirichlet fundamental polyhedron $D_z(\Gamma)$. This is because the new polyhedron $D_\eta \subset H^2_\mathbb{F}$ has isometrically the same (tessellations of) neighborhoods of its edges (= intersections of its sides) as $D_z(\Gamma)$ had. This implies that the polyhedra $\gamma'(D_\eta)$, $\gamma' \in \Gamma_\eta$, form a tessellation of $H^2_\mathbb{F}$ (with non-overlapping interiors). Hence the deformed group $\Gamma_\eta \subset \operatorname{Isom} H^2_\mathbb{F}$ is a discrete group, and $D_\eta$ is its fundamental polyhedron bounded by bisectors.

$\Gamma$-equivariantly extending our elementary bending homeomorphism between two closed polyhedra, $\varphi_{\eta,\zeta} : \overline{D_z(\Gamma)} \to \overline{D_\eta}$, we obtain a desired $\Gamma$-equivariant homeomorphism $F_\eta : \overline{H^2_\mathbb{F}} \backslash \Lambda(\Gamma) \to \overline{H^2_\mathbb{F}} \backslash \Lambda(\Gamma_\eta)$. Its extension by continuity to the limit (real) circle $\Lambda(\Gamma)$ coincides with the canonical equivariant homeomorphism $f_\chi : \Lambda(\Gamma) \to \Lambda(\Gamma_\eta)$ given by the isomorphism $\chi : \Gamma \to \Gamma_\eta$, see [2]. Hence we have a $\Gamma$-equivariant self-homeomorphism of the whole space $\overline{H^2_\mathbb{F}}$, which we denote as before by $F_\eta$.

We claim that for a small $\eta_0 > 0$, the family $\{F^*_\eta\}$ of representations $F^*_\eta : \Gamma \to \Gamma_\eta = F_\eta \Gamma F_\eta^{-1}$, $\eta \in (-\eta_0, \eta_0)$, defines a nontrivial curve $\mathcal{B} : (-\eta_0, \eta_0) \to \mathcal{R}_0(\Gamma)$ in the representation variety $\mathcal{R}_0(\Gamma)$, which covers a curve in the Teichmüller space $\mathcal{T}(\Gamma)$ represented by conjugacy classes of nontrivial quasifuchsian representations $[\mathcal{B}(\eta)] = [F^*_\eta]$. We call the constructed deformation $\mathcal{B}$ a *bending deformation* of a given discrete group $\Gamma \subset PO(2,1) \subset \operatorname{Isom} H^2_\mathbb{F}$ along a (small) bending geodesic. Since our parameter $\eta$ represents elements in the compact group $\mathcal{U}_\mathbb{F}$ of unit $\mathbb{F}$-numbers with real dimension $k_\mathbb{F} = \dim \mathcal{U}_\mathbb{F}$ ($k_\mathbb{F}$ is 1, 3 and 7 respectively for complex, quaternionic and octonionic cases), our deformation curve is $k_\mathbb{F}$-dimensional. So we have proved:

**Theorem 6.2** *Let $\Gamma \subset PO(2,1) \subset \operatorname{Isom} H^2_\mathbb{F}$ be a given discrete group uniformizing a Riemann surface $S = H^2_\mathbb{R}/\Gamma$ with a small simple closed geodesic $\alpha \subset S$. Then for a sufficiently small $\eta_0 > 0$, there is a bending deformation $\mathcal{B}_\alpha : B^{k_\mathbb{F}}(0, \eta_0) \to \mathcal{R}_0(\Gamma)$ of the group $\Gamma$ along $\alpha$, $\mathcal{B}_\alpha(\eta) = \rho_\eta = F^*_\eta$, induced by* $\Gamma$-*equivariant homeomorphisms* $F_\eta : \overline{H_\mathbb{F}^2} \to \overline{H_\mathbb{F}^2}$, $\eta \in B^{k_\mathbb{F}}(0, \eta_0)$, *with* $k_\mathbb{F} = 1, 3$ *and* $7$ *for* $\mathbb{F} = \mathbb{C}, \mathbb{H}$ *and* $\mathbb{O}$, *respectively.*

Noticing that we may deform $\Gamma$ in $PO(2,1)$ to have $3g-3$ small disjoint closed geodesics on $S_g = H_\mathbb{R}^2/\Gamma$ (Prop. 6.1), we see from (10) that bendings along disjoint closed geodesics are independent. This and the fact that the Teichmüller space $\mathcal{T}(S_g)$ is an open cell in $\mathbb{C}^{3g-3}$ imply:

**Corollary 6.3** *Let* $S_g = H_\mathbb{R}^2/\Gamma$ *be a closed totally real geodesic surface of genus* $g > 1$ *in a given* $\mathbb{F}$-*hyperbolic manifold* $M = H_\mathbb{F}^2/\Gamma$, $\Gamma \subset PO(2,1) \subset$ Isom $H_\mathbb{F}^2$. *Then there is a real analytic embedding* $\mathcal{B} : B^{(6+3k_\mathbb{F})(g-1)} \hookrightarrow \mathcal{T}(M)$ *of a real* $(6+3k_\mathbb{F})(g-1)$-*ball into the Teichmüller space of* $M$, *defined by deformations along disjoint closed geodesics in* $M$, *with* $k_\mathbb{F} = 1, 3$ *and* $7$ *for* $\mathbb{F} = \mathbb{C}, \mathbb{H}$ *and* $\mathbb{O}$, *respectively.*

**Proof:** First, we observe that our construction works not only in the case of $\mathbb{F}$-hyperbolic manifolds $M$ homotopy equivalent to totally real geodesic surface but for manifolds $M$ homotopy equivalent to a 2-surface which has a totally real geodesic piece bounded by closed geodesics. In such a case, the holonomy group of that piece is a non-elementary discrete subgroup in Isom $H_\mathbb{F}^2$ preserving a totally real geodesic plane $H_\mathbb{R}^2 \subset H_\mathbb{F}^2$. Then we see from (10) that bendings along disjoint closed geodesics are independent. In addition to the above construction, we have also to show that our bending deformation is not trivial, and $[\mathcal{B}(\eta)] \neq [\mathcal{B}(\eta')]$ for any $\eta \neq \eta'$.

The non-triviality of this deformation follows directly from (10), cf. [1]. Namely, the restrictions $\rho_\eta|_{\Gamma_1}$ of bending representations to a non-elementary subgroup $\Gamma_1 \subset \Gamma$ (in general, to a "real" subgroup $\Gamma_r \subset \Gamma$ corresponding to a totally real geodesic piece in the homotopy equivalent surface $S \simeq M$) are identical. So if the deformation $\mathcal{B}$ were trivial then it would be conjugation of the group $\Gamma$ by $\mathbb{F}$-isometries that commute with the non-trivial real subgroup $\Gamma_r \subset \Gamma$ and pointwise fix the totally real geodesic plane $H_\mathbb{R}^2$. This contradicts to the fact that the limit set of any deformed group $\Gamma_\eta$, $\eta \neq 0$, does not belong to the real circle containing the Cantor limit set $\Lambda(\Gamma_r)$.

The injectivity of the map $\mathcal{B}$ we obtain by using the quaternionic and octonionic angular invariants $\mathbb{A}_\mathbb{F}(x)$ for triples $x = (x^0, x^1, x^2)$ of points in $\partial H_\mathbb{F}^2 = N \cup \{\infty\}$, see Definitions 3.1 and 5.1. Namely, let $\gamma_2 \in \Gamma \backslash \Gamma_1$ be a generator of the group $\Gamma$ in (9) whose fixed point $x^2 \in \Lambda(\Gamma)$ lies in $\mathbb{R}_+ \times \{0\} \subset N$, and $x^2_\eta \in \Lambda(\Gamma_\eta)$ the corresponding fixed point of the element $\chi_\eta(\gamma_2) \in \Gamma_\eta$ under the free-product isomorphism $\chi_\eta \colon \Gamma \to \Gamma_\eta$. Due to our construction, one can see that the orbit $\gamma(x_\eta^2)$, $\gamma \in \Gamma_P$, under the loxodromic (dilation) subgroup $\Gamma_P \subset \Gamma \cap \Gamma_\eta$ approximates the origin along a ray $(0, \infty)$ which has a non-zero angle $\eta$ with the ray $\mathbb{R}_- \times \{0\} \subset N$. The latter ray also contains an orbit $\gamma(x^1)$, $\gamma \in \Gamma_P$, of a limit point $x^1$ of $\Gamma_1$ which approximates the origin from the other side. Taking triples $x = (x^1, 0, x^2)$ and $x_\eta = (x^1, 0, x_\eta^2)$ of points which lie correspondingly in the limit sets $\Lambda(\Gamma)$ and $\Lambda(\Gamma_\eta)$, we have that $\mathbb{A}_\mathbb{F}(x) = 0$ and $\mathbb{A}_\mathbb{F}(x_\eta) \neq 0, \pm\pi/2$. Therefore these limit sets (topological circles) cannot be equivalent under any $\mathbb{F}$-hyperbolic isometry because of different angular invariants (another view on nontriviality of deformations).

Similarly, for two different values $\eta$ and $\eta'$, we have triples $x_\eta$ and $x_{\eta'}$ with different (non-trivial) Cartan angular invariants $\mathbb{A}_\mathbb{F}(x_\eta) \neq \mathbb{A}_\mathbb{F}(x_{\eta'})$. Hence two topological circles $\Lambda(\Gamma_\eta)$ and $\Lambda(\Gamma_{\eta'})$ are not $\mathbb{F}$-hyperbolically equivalent. This completes the proof of injectivity. ∎

## 6.3 Deformations of non-surface groups

As for quaternionic case, we note that the octonionic angular invariant $\mathbb{A}_\mathbb{O}(x)$ cannot distinguish different shapes of $k$-spheres inside the sphere at infinity of an $\mathbb{O}$-line, in particular, a $k$-dimensional linear subspace $L$ and the result $L'$ of its bending along a hyperplane, see next Section. In both cases, for any triple $x \in L$ and the corresponding triple $x' \in L'$, we have $\mathbb{A}_\mathbb{O}(x) = \mathbb{A}_\mathbb{O}(x') = \pi/2$. Hence $\mathbb{A}_\mathbb{O}(x)$ cannot distinguish a Fuchsian group $\Gamma \subset \text{Isom } H_\mathbb{O}^1$ from its non-trivial quasi-Fuchsian deformation $\rho(\Gamma) \subset \text{Isom } H_\mathbb{O}^1$. In particular, we may construct such non-trivial quasi-Fuchsian deformations as follows.

**Proposition 6.4** *Let $\Gamma$ be a lattice in $PSp(1,1) \cong \text{Isom}_+ H_\mathbb{R}^4$ such that the quotient $H_\mathbb{H}^1/\Gamma$ has a closed totally geodesic 3-manifold. Then there exists a continuous non-trivial quasi-Fuchsian deformation*

$$R \colon (-\varepsilon, \varepsilon) \to \text{Hom}(\Gamma, F_4^{-20}), \quad R(\varepsilon) \colon \Gamma \to \text{Isom } H_\mathbb{O}^1 \subset \text{Isom } H_\mathbb{O}^2.$$

**Proof:** Since $H_\mathbb{O}^1 \cong H_\mathbb{R}^8$, we may assume that $\Gamma \subset \text{Isom}_+ H_\mathbb{R}^4$ preserves a totally geodesic subspace $H_\mathbb{R}^4 \subset H_\mathbb{R}^8$ and has a bending hyperplane $P \subset H_\mathbb{R}^4$. Now we take a stabilizer $Spin_P \subset Spin(7)$ of this 3-dimensional plane $P$. It is isomorphic to $SO(5)$ and acts transitively on the unit 4-sphere in the orthogonal complement (in $H_\mathbb{O}^1$) to $P$. Elements $U_\eta \in Spin_P$ bend the totally geodesic subspace $H_\mathbb{R}^4$ (inside $H_\mathbb{O}^1$) along the plane $P$. Therefore, in both cases of a free amalgamated product or $HNN$-group $\Gamma$, see (9), we have non-trivial bending representations $\rho_\eta : \Gamma \to \Gamma_\eta \subset F_4^{-20}$, where deformed groups $\Gamma_\eta$ still preserve the $\mathbb{O}$-line $H_\mathbb{O}^1$ and are given in (10). ∎

In contrast to such a flexibility, Corlette's [12] superrigidity implies (cf.[20]):

**Theorem 6.5** *Let $\Gamma$ be a uniform lattice in $PSp(2,1)$, and $G$ be either $PSp(m,1)$, $m \geq 2$, or $F_4^{-20}$. Then any discrete representation $\rho : \Gamma \to G$ is conjugate to the natural inclusion $\Gamma \subset G$. In particular, $\rho(\Gamma)$ preserves a quaternionic 2-subspace $H_\mathbb{H}^2$ (either in $H_\mathbb{H}^m$ or $H_\mathbb{O}^2$).*

**Proof:** Let $X$ be the smallest totally geodesic subspace in $H_\mathbb{O}^2$ whose infinity contains the limit set $\Lambda$ of $\rho(\Gamma)$. Such a subspace $X$ can be described by:

**Lemma 6.6** *Let $X$ be a totally geodesic subspace in the hyperbolic Cayley plane $H_\mathbb{O}^2$. Then $X$ is isometric either to a subspace $H_\mathbb{R}^k \subset H_\mathbb{O}^1$, $1 \leq k \leq 8$, or to $H_\mathbb{F}^2$ where $\mathbb{F}$ is $\mathbb{R}, \mathbb{C}, \mathbb{H}$ or $\mathbb{O}$. Furthermore, the stabilizer $G_X$ of $X$ in $F_4^{-20}$ is $\mathrm{Aut}(X) \times K'$ where $K' \subset Spin(9)$ pointwise fixes $X$.*

Clearly $\rho(\Gamma)$ preserves $X$ due to the invariance of the limit set $\Lambda$, so we may think that $\rho(\Gamma) \subset \mathrm{Aut}(X) \subset F_4^{-20}$. Since $\rho(\Gamma)$ is Zariski dense in $\mathrm{Aut}(X)$, the Corlette's superrigidity implies that $X \neq H_\mathbb{O}^2$. Since $\rho$ is discrete, faithful and Zariski dense in $\mathrm{Aut}(X)$, the Corlette's superrigidity implies that $\rho$ extends to a homomorphism $\rho : \mathrm{Aut}(H_\mathbb{H}^2) \to \mathrm{Aut}(X)$. However $\mathrm{Aut}(H_\mathbb{H}^2)$ is semisimple with trivial center, so the obtained homomorphism is an isomorphism. This is possible only when $X = H_\mathbb{H}^2$. ∎

## 7 Representations of surface groups

Let $S$ be a closed Riemann surface of genus $g > 1$ and $\rho : \pi_1(S) \to PSp(n,1)$ be a representation. Let $X_\rho$ be the associated flat $H_\mathbb{H}^n$ bundle over $S$ with holonomy $\rho$, i.e., it is the quotient of $\tilde{S} \times H_\mathbb{H}^n$ by the diagonal action of $\pi_1(S)$ via $\rho$. A section of $X_\rho$ defines a $\pi_1(S)$ equivariant map $\tilde{f}$ from $\tilde{S}$ to $H_\mathbb{H}^n$. Define the straight 2-cochain $\tilde{f}^*\omega$ by

$$\tilde{f}^*\omega\langle x_1, x_2, x_3\rangle = \int_{\triangle_{\tilde{f}(x)}} \omega = 4\pi\tau_\mathbb{H}(\tilde{f}(x)),$$

where $\tilde{f}(x) = (\tilde{f}(x_1), \tilde{f}(x_2), \tilde{f}(x_3))$. The cochain $\tilde{f}^*\omega$ is $\pi_1(S)$-equivariant, so descends to $f^*\omega$ on $S$. For the orthogonal projection $\Pi$ onto an $\mathbb{H}$-line containing points $x_1$ and $x_2$, we have

$$\left| \int_{\triangle(x_1,x_2,x_3)} \omega \right| = |\operatorname{Area}(\triangle(x_1, x_2, \Pi(x_3))) \cdot \operatorname{Area}(S^2)| \leq 4\pi^2,$$

Therefore, $f^*\omega$ represents a bounded cohomology class on $S$ with sup norm at most $4\pi^2$. Note that the equality holds only if the points $x_1, x_2, x_3$ lie on the boundary of an $\mathbb{H}$-line.

Now we define an invariant $c(\rho)$ as the character of the representation $\rho$, which is the evaluation of $f^*\omega$ on the fundamental cycle $[S]$.

## 7.1 Hybrid representations

Studying which integers may be the values of the complex Toledo invariant, [17] constructed a class of specific complex surfaces that are disc bundles over a closed surface $S_g$ of genus $g > 1$ and showed that every even integer number $\tau$, $2 - 2g \leq \tau \leq 2g - 2$, is realized as the complex Toledo invariant of such a complex surface $M$ (of its holonomy representation $\rho_\mathbb{C}: \pi_1 M \to PU(2, 1)$).

Actually it was constructed a series of discrete groups $\Gamma \subset PU(2, 1)$ such that each quotient manifold $M = H_\mathbb{C}^2/\Gamma$ is homotopy equivalent to a closed piecewise hyperbolic surface $S$ of genus $g > 1$ with collection $\Sigma_0$ of simple closed geodesics such that on each component of $S \backslash \Sigma_0$ the metric either has constant curvature $-1$ (such components form a subsurface $\Sigma_1 \subset S \backslash \Sigma_0$) or constant curvature $-1/4$ (such components form a subsurface $\Sigma_2 \subset S \backslash \Sigma_0$). Therefore $S = \Sigma_1 \cup \Sigma_0 \cup \Sigma_2$ where $\Sigma_1$ is holomorphically geodesic and $\Sigma_2$ is totally real. We may represent the fundamental groups of these components correspondingly as discrete (Fuchsian) subgroups $\Gamma_1 \subset PU(1, 1)$ and $\Gamma_2 \subset PO(2, 1)$, so that the group $\Gamma$ can be obtained from them by a sequence of free amalgamated products and HNN-extensions. In terms of complex surfaces $M_i = H_\mathbb{C}^2/\Gamma_i$, this construction represents gluing of them and is called *hybriding*. In this way, one can obtain a discrete embedding $\rho_\mathbb{C}: \Gamma \hookrightarrow PU(2, 1)$ which is the holonomy representation of a complex surface $M$ and whose complex Toledo invariant is the given even number $\tau$.

From the previous section we know that for a given (quaternionic) representation $\rho: \pi_1(S) \to PSp(n, 1)$, its character $c(\rho) = 0$ if the holonomy group $\rho(\pi_1(S))$ leaves invariant a totally real subspace in $H_\mathbb{H}^n$, and $|c(\rho)|$ has its maximum when $\rho(\pi_1(S))$ has an invariant quaternionic line in $H_\mathbb{H}^n$. Now we observe that each real line $(\mathbb{R}, 0) \subset \mathbb{F}^{n-1} \times \operatorname{Im} \mathbb{F}$ passing through the origin lies at infinity of a totally real plane $H^2_\mathbb{R} \subset H^n_\mathbb{F}$, and each one dimensional line in $(0, \operatorname{Im} \mathbb{F})$ passing through the origin lies at infinity of a complex geodesic in $H^n_\mathbb{H}$. That is why the above (complex hyperbolic) hybriding construction perfectly work in the case of hyperbolic quaternionic (and octonionic) spaces. Directly following that hybriding construction, we immediately have a discrete representation $\rho : \pi_1(S) \to PSp(n, 1)$, whose characteristic value $c(\rho)$ is between two extremes, namely:

$$|c(\rho)| = 8\pi^2 |\chi(\Sigma_1)|. \tag{12}$$

## 7.2 Mixture of hybriding and bending

Now, continuing the hybriding construction as in the above section, we take a closed geodesic $\alpha$ in the "real" part $\Sigma_2$ (with curvature $-1/4$) of our piecewise hyperbolic surface $S \subset H^n_\mathbb{F}/\Gamma$, $\pi_1(S) \cong \Gamma$, and bend the whole group $\Gamma$ along $\alpha$, by a small enough angle $\eta$, to obtain a new discrete faithful representation $\rho_\eta : \Gamma \to \Gamma_\eta \subset \operatorname{Isom} H^n_\mathbb{F}$. Clearly, the character $c(\rho_\eta)$ of the obtained representation is the same as for the previous one in (12), that is for the inclusion $\Gamma \subset \operatorname{Isom} H^n_\mathbb{F}$, but those representations are not conjugate each other any more.

The above situation is just one of possible applications of bending deformations in the case of $\mathbb{F}$-hyperbolic structures. One can apply the proof of Corollary 6.3 to a general situation of bending deformations of a $\mathbb{F}$-hyperbolic manifold $M = H^2_\mathbb{F}/\Gamma$ whose holonomy group $\Gamma \subset \operatorname{Isom} H^2_\mathbb{F}$ has a non-elementary subgroup $\Gamma_r$ preserving a totally real geodesic plane $H^2_\mathbb{R}$ (where some convex subdomain is precisely $\Gamma_r$-invariant in $\Gamma$). In other words, such a manifold $M$ has an embedded totally real geodesic surface with geodesic boundary. So we immediately have:

**Corollary 7.1** *Let $M = H^2_\mathbb{F}/\Gamma$ be an $\mathbb{F}$-hyperbolic manifold with embedded totally real geodesic surface $S_r \subset M$ with geodesic boundary. Then the bending deformation $\mathcal{B} : (-\eta, \eta) \to \mathcal{T}(M)$ of $M$ along a simple closed geodesic $\alpha \subset S_r$ is a real analytic embedding provided that the limit set $\Lambda(\Gamma)$ of the holonomy group $\Gamma$ does not belong to the $\Gamma$-orbit of a real circle $S^1_\mathbb{R}$ and a complex chain $S^1_\mathbb{C}$, where the latter is the infinity of the complex geodesic containing a lift $\tilde{\alpha} \subset H^2_\mathbb{F}$ of the closed geodesic $\alpha$, and the former one contains the limit set of the holonomy group $\Gamma_r \subset \Gamma$ of the totally real geodesic surface $S_r$.*

Dept of Math., University of Oklahoma, Norman, OK 73019-0351, USA
Dept of Math., Seoul National University, Seoul, Korea
e-mail: apanasov@ou.edu, inkang@math.snu.ac.kr